\pdfoutput=1
\documentclass[nohyperdvips, review]{siamonline190516}
\usepackage{layouts}

\usepackage{lipsum}
\usepackage{amsfonts}
\usepackage[mathscr]{euscript} % <-- SB preferred script \mathscr
\usepackage{dsfont} % <-- blackboard bold 1 (one)
\usepackage{graphicx}
\usepackage{epstopdf}
\usepackage[noend]{algorithmic}

\usepackage{enumitem}
\usepackage{comment}

\ifpdf
  \DeclareGraphicsExtensions{.eps,.pdf,.png,.jpg}
\else
  \DeclareGraphicsExtensions{.eps}
\fi

\usepackage{cite}

\usepackage{inconsolata}

\usepackage[multiple]{footmisc}

\newcommand{\Real}{\mathbb{R}}

\newcommand{\Cinfty}{\ensuremath{{C}^\infty}}
\newcommand{\Cr}{\ensuremath{{C}^r}}
\newcommand{\Co}{\ensuremath{{C}^1}}
\newcommand{\ECr}{\ensuremath{{E}\Cr}}
\newcommand{\PCr}{\ensuremath{{P}\Cr}}

\newcommand{\Frechet}{Fr{\'{e}}chet}
\newcommand{\naive}{na{\"{i}}ve}

\newcommand{\ones}{\mathds{1}}

\newcommand{\figref}[1]{\cref{fig:#1}}
\newcommand{\secref}[1]{\cref{sec:#1}}
\newcommand{\algref}[1]{\cref{alg:#1}}

\newcommand{\R}{\mathbb{R}}
\newcommand{\N}{\mathbb{N}}
\newcommand{\Ord}[1]{\ensuremath{O\paren{#1}}}
\newcommand{\perm}[1]{\ensuremath{\sigma(\set{0,\dots,#1})}}

\newcommand{\F}{\mathscr{F}} % SB flow domain
\newcommand{\K}{\mathscr{K}} % SB kernel of Dh(rho)
\newcommand{\Li}{\mathscr{L}} % SB lineality space of conical subdivision Sigma'
\newcommand{\J}{J} % SB index set
\newcommand{\Sn}{{S}_n} % SB symmetric group of order m
\newcommand{\Bn}{{B}_n} % SB corners of n-dimensional hypercube 
\newcommand{\Bd}{{B}_d} % SB corners of d-dimensional hypercube
\newcommand{\tr}{\top} % SB I prefer this transpose -- thoughts?
\newcommand{\pinv}{\dagger}

\newcommand{\into}{\rightarrow}
\newcommand{\assign}{\leftarrow}
\newcommand{\goesto}{\rightarrow}
\newcommand{\Int}{\operatorname{Int}}
\newcommand{\conv}{\operatorname{conv}}
\newcommand{\cone}{\operatorname{cone}}
\newcommand{\spn}{\operatorname{span}}

\newcommand{\aff}{\operatorname{aff}}
\renewcommand{\dim}{\operatorname{dim}}
\newcommand{\rank}{\operatorname{rank}}
\newcommand{\sgn}{\operatorname{sign}}
\newcommand{\sm}{\setminus}
\newcommand{\abs}[1]{\left| #1 \right|}
\newcommand{\car}[1]{\left| #1 \right|}
\newcommand{\card}[1]{\#(#1)}
\newcommand{\norm}[1]{\left\| #1 \right\|}
\newcommand{\paren}[1]{\left( #1 \right)}
\newcommand{\brak}[1]{\left[ #1 \right]}

\newcommand{\mat}[2]{\brak{\begin{array}{#1} #2 \end{array}}}
\newcommand{\pw}[1]{\begin{cases} #1 \end{cases}}
\newcommand{\td}[1]{\widetilde{#1}}
\newcommand{\samp}[1]{\td{#1}}
\newcommand{\set}[1]{\left\{#1\right\}}
\newcommand{\eqn}[1]{\begin{equation*}\begin{aligned} #1 \end{aligned}\end{equation*}}
\newcommand{\eqnn}[1]{\begin{equation}\begin{aligned} #1 \end{aligned}\end{equation}}

\newsiamremark{remark}{Remark}
\newsiamremark{hypothesis}{Hypothesis}
\crefname{hypothesis}{Hypothesis}{Hypotheses}
\newsiamthm{claim}{Claim}
\newsiamremark{example}{Example}

\newcommand{\thetitle}{Representing and computing the B-derivative of an {\ECr} vector field's {\PCr} flow}
\newcommand{\theheadertitle}{B-derivative of an $\protect\ECr$ vector field's $\protect\PCr$ flow}
\newcommand{\theauthors}{G. Council, S. Revzen, S. A. Burden}

\headers{\theheadertitle}{\theauthors}

\title{\thetitle\thanks{Submitted to the editors DATE.
\funding{This material is based upon work supported by the U.~S.~Army Research Laboratory and the U.~S.~Army Research Office under contract/grant number W911NF-16-1-0158, ARO W911NF-14-1-0573, and ARO MURI W911NF-17-1-0306. It is also supported by U.~S.~National Science Foundation Cyber-Physical Systems Award~\#1836819 and National Robotics Initiative Award~\#1924303.}}}

\author{%
George Council\thanks{Department of Mechanical Engineering, Carnegie Mellon University, Pittsburgh, PA, USA (\email{gcouncil@andrew.cmu.edu}).}
\and Shai Revzen\thanks{Department of Electrical Engineering and Computer Science, University of Michigan, Ann Arbor, MI, USA (\email{shrevzen@umich.edu}).}
\and Samuel A. Burden\thanks{Department of Electrical \& Computer Engineering, University of Washington, Seattle, WA, USA (\email{sburden@uw.edu}, \url{http://faculty.uw.edu/sburden}).}
}

\usepackage{amsopn}

\graphicspath{{figures/}}
\usepackage{pdfpages}
\providecommand*{\backrefsetup}[1]{}

\nolinenumbers
\begin{document}
\maketitle

% REQUIRED
\begin{abstract}
  This paper concerns 
  first-order approximation of the piecewise-differentiable flow generated by a class of nonsmooth vector fields.
  Specifically, we represent and compute the Bouligand (or B-)derivative of the piecewise-{\Cr} flow generated by an event-selected {\Cr} vector field.
  Our results are remarkably efficient:
  although there are factorially many ``pieces'' of the desired derivative,
  we provide an algorithm that evaluates its action on a given tangent vector using polynomial time and space,
  and verify the algorithm's correctness by deriving a representation for the B-derivative that requires ``only'' exponential time and space to construct.
  We apply our methods in two classes of illustrative examples: piecewise-constant vector fields and mechanical systems subject to unilateral constraints.
\end{abstract}

% REQUIRED
\begin{keywords}
nonsmooth dynamical system, 
differential equation with discontinuous right-hand side, 
first-order approximation,
Bouligand derivative,
saltation matrix,
\end{keywords}

% REQUIRED
\begin{AMS}
34A36, % Discontinuous ordinary differential equations
%34A37, % Differential equations with impulses
%34A38, % Hybrid systems
65D30, % Numerical integration
65D99, % None of the above, but in this section:  Numerical approximation and computational geometry (primarily algorithms)
%68T40, % Robotics
70E99, % Robot dynamics and control
\end{AMS}

\section{Introduction}
\label{sec:intro}
First-order approximations -- i.e.\ derivatives -- are a foundational tool for analysis and synthesis in smooth dynamical and control systems.
For instance, derivatives play a crucial r{\^{o}}le in: % theoretical and computational aspects of:
%\begin{itemize}
%  \item stability via spectral~\cite[Ch.~8.3]{Hespanha2009-nf} or Lyapunov~\cite[Ch.~5]{Sastry1999-ei} methods;
%  \item controllability via linearization~\cite[Ch.~8.7]{Hespanha2009-nf} or Frobenius/Chow~\cite[Ch.~8/Ch.~11]{Sastry1999-ei} techniques;
%  \item optimality via stationarity~\cite[Ch.~1]{Bertsekas1999-an} or Pontryagin~\cite[Ch.~1]{Pontryagin1962-rn} principles;
%  \item identifiability via adaptation~\cite[Ch.~2]{Sastry1989-zr} or Expectation-Maximization~\cite[Ch.~10]{Ljung1999-ee}.
%\end{itemize}
stability via spectral~\cite[Ch.~8.3]{Hespanha2009-nf} or Lyapunov~\cite[Ch.~5]{Sastry1999-ei} methods;
controllability via linearization~\cite[Ch.~8.7]{Hespanha2009-nf} or Frobenius/Chow~\cite[Ch.~8/Ch.~11]{Sastry1999-ei} techniques;
optimality via stationarity~\cite[Ch.~1]{Bertsekas1999-an} or Pontryagin~\cite[Ch.~1]{Pontryagin1962-rn} principles;
identifiability via adaptation~\cite[Ch.~2]{Sastry1989-zr} or Expectation-Maximization~\cite[Ch.~10]{Ljung1999-ee} methods.
These tools all depend on the existence of a computationally-amenable representation for the first-order approximation of smooth system dynamics -- namely, the \emph{\Frechet} (or F-)derivative of the system's smooth flow~\cite[Ch.~5.6]{Polak1997-xd}, which derivative is a continuous linear function of tangent vectors.%
\footnote{
We emphasize both properties of the {\Frechet} derivative (continuity \emph{and} linearity) since the generalized derivative we consider in what follows will retain one (continuity) while relaxing the other (piecewise-linearity).
%A linear function between finite-dimensional vector spaces is always continuous, but a linear function between Banach spaces is continuous if and only if its induced norm is finite~\cite[Ch.~2.9]{Luenberger1969-ha}.
}

By definition, nonsmooth systems do not generally enjoy existence (let alone computational amenability) of first-order approximations.
Restricting to the class of (so-called~\cite[Def.~1,~2]{Burden2016-bb}) \emph{event-selected \Cr} (\ECr) vector fields that (i) are smooth except along a finite number of surfaces of discontinuity and 
(ii) preclude \emph{sliding}~\cite{Utkin1977-rm, Jeffrey2014-nt} or \emph{branching}~\cite[Def.~3.11]{Simic2005-fv} through a transversality condition, 
we obtain flows that are piecewise-differentiable~\cite[Thm.~4]{Burden2016-bb} (specifically, \emph{piecewise}-\Cr ({\PCr})~\cite[Ch.~4.1]{Scholtes2012-la}).
By virtue of their piecewise-differentiability,
these flows admit a first-order approximation, termed the \emph{Bouligand} (or B-)derivative, which derivative is a continuous \emph{piecewise}-linear function of tangent vectors~\cite[Ch.~3,~4]{Scholtes2012-la}.
This paper is concerned with the efficient representation and computation of this piecewise-linear first-order approximation.

Our contributions are twofold:
(i) we construct a representation for the B-derivative of the {\PCr} flow generated by an {\ECr} vector field;
(ii) we derive an algorithm that evaluates the B-derivative on a given tangent vector.
Although there are factorially many ``pieces'' of the derivative,
% our methods are remarkably efficient:
we 
(i) represent it using exponential time and space
and
(ii) compute it using polynomial time and space. 
%; as a side-effect of our analysis, we show that these complexities cannot be improved upon in general.
In an effort to make our results as accessible and useful as possible, 
we 
provide a concise summary of the algorithm in~\secref{alg} 
and
apply our methods in~\secref{app} 
%to 
%piecewise-constant vector fields %in~\secref{app:pwc}
%and 
%mechanical systems subject to unilateral constraints. %in~\secref{app:mech}.
\emph{before} 
rehearsing the technical background in~\secref{bg} needed to derive the representation in~\secref{rep} and verify the algorithm's correctness in~\secref{comp}.

We emphasize that our methods are most useful when there are more than two surfaces of discontinuity, as representation and computation of first-order approximations in the $1$- and $2$-surface cases have been investigated extensively~\cite{Aizerman1958-ih, Ivanov1998-ff, Hiskens2000-ps, Dieci2011-ps, Bernardo2008piecewise, Bizzarri2013-jh}, and these cases do not benefit from the complexity savings touted above.
Previously, we established existence of the piecewise-linear first-order approximation of the flow~\cite[Rem.~1]{Burden2016-bb} and provided an inefficient scheme to evaluate each of its ``pieces''~\cite[Sec.~7]{Burden2016-bb} in the presence of an arbitrary number of surfaces of discontinuity.
To the best of our knowledge, the present paper contains the first representation for the B-derivative of the {\PCr} flow of a general {\ECr} vector field and polynomial-time algorithm to compute it.

\section{Algorithm}
\label{sec:alg}
The goal of this paper is to obtain an algorithm that efficiently computes the derivative of a class of nonsmooth flows.
This computational task and our solution are easy to describe, yet verifying the algorithm's correctness requires significant technical overhead.
Thus, the remainder of this section will be devoted to specifying the algorithm and the problem it solves using minimal notation and terminology. 
Subsequent sections will provide technical details -- which may be of interest in their own right -- that prove the algorithm is correct.

%\subsection{Problem statement}

Given vector field $F:\R^d\into T\R^d$ and trajectory $x:[0,\infty)\into\R^d$ satisfying%
\footnote{%
In this section, we will denote time dependence using subscripts rather than parentheses.
}
\eqnn{\label{eq:alg:trj}
\forall t \geq 0 : x_t = \int_0^t F(x_\tau) \,d\tau,
}
our goal is to approximate how $x_t$ varies with respect to $x_0$ to first order for a given $t > 0$.
Formally, with $\phi:[0,\infty)\times\R^d\into\R^d$ denoting the \emph{flow} of $F$ satisfying
\eqnn{\label{eq:alg:flow}
\forall t \geq 0, x_0\in\R^d : \phi_t(x_0) = \int_0^t F\paren{\phi_\tau(x_0)} \,d\tau,
}
our goal is to evaluate the directional derivative 
$D \phi_t\paren{x_0; \delta x_0}$ 
given $t > 0$, $\delta x_0\in T_{x_0} \R^d$:
\eqnn{\label{eq:dder}
\forall t > 0, \delta x_0\in T_{x_0} \R^d: D \phi_t\paren{x_0; \delta x_0} = \lim_{\alpha\goesto 0^+} \frac{1}{\alpha}\paren{\phi_t(x_0 + \alpha\ \delta x_0) - \phi_t(x_0)}.
}
Specifically, we seek to evaluate this derivative for vector fields that are smooth everywhere except a finite collection of surfaces where they are allowed to be discontinuous.
We will first recall how to obtain the derivative in the presence of zero 
(\secref{alg:0})
and
one (\secref{alg:1})
surfaces of discontinuity 
before presenting our algorithm, which is applicable in the presence of an arbitrary number of surfaces of discontinuity (\secref{alg:n}).

\subsection{Continuously-differentiable vector field}
\label{sec:alg:0}
If $F$ is continuously differentiable on the trajectory $x$, 
the derivative 
$\delta x_t = D \phi_t\paren{x_0; \delta x_0}$ 
satisfies the
linear time-varying
\emph{variational equation}~\cite[Appendix~B]{Parker1989-df} %~\cite[Sec.~15.2, Thm.~1]{Hirsch1974-lt}
\eqnn{\label{eq:alg:var}
\forall t \geq 0 : \delta x_t = \int_0^t D F(x_\tau)\cdot\delta x_\tau \,d\tau,
}
whence 
$\delta x_t = D \phi_t\paren{x_0; \delta x_0}$ 
can be approximated to any desired precision in polynomial time by applying numerical simulation algorithms~\cite[Ch.~4]{Parker1989-df} to~\eqref{eq:alg:trj},~\eqref{eq:alg:var}.

\subsection{Single surface of discontinuity}
\label{sec:alg:1}
If $F$ is continuously differentiable everywhere except a smooth codimension-1 submanifold $H\subset\R^d$ that intersects the trajectory $x$ transversally at only one point $x_{s}$, $s\in(0,t)$,
the continuous-time equation~\eqref{eq:alg:var} is augmented by the discrete-time update~\cite[Eqn.~(58)]{Aizerman1958-ih}, 
\eqnn{\label{eq:alg:salt}
\delta x_s^+ = \paren{ I_d + \frac{\paren{F^+ - F^-}\cdot \eta^\tr}{\eta^\tr\cdot F^-} } \cdot \delta x_s^- = M\cdot\delta x_s^-,
}
where 
$\delta x_s^\pm = \lim_{\tau\goesto s^\pm} \delta x_\tau$
and
$F^\pm = \lim_{\tau\goesto s^\pm} F\paren{x_\tau}$
denote the limiting values of 
$\delta x_\tau$ 
and
$F\paren{x_\tau}$
at $s$ 
from the right ($+$) and left ($-$)
and
$\eta\in\R^d$ 
is any vector orthogonal to surface $H$ at $x_s$;
$M\in\R^{d\times d}$ is termed the \emph{saltation matrix} \cite[Eqn.~(2.76)]{Bernardo2008piecewise}, \cite[Eqn.~(7.65)]{Leine2013dynamics}.
Overall, the desired derivative is
\eqnn{
D \phi_t(x_0; \delta x_0) = D \phi_{t-s}(x_s) \cdot M \cdot D \phi_s(x_0) \cdot \delta x_0,
}
where $D \phi_{t-s}(x_s), D \phi_s(x_0)\in\R^{d\times d}$ can be approximated by simulating~\eqref{eq:alg:flow},~\eqref{eq:alg:var} since the flow is smooth away from time $s$.
Computing the saltation matrix $M$ requires $\Ord{d^2}$ time and space,
but evaluating its action on $\delta x_s^-$ in~\eqref{eq:alg:salt} requires only $\Ord{d}$ time and space.
%
%\SB{is the following worth including?  I feel like it is easy enough for the interested reader to verify, so opted to exclude it initially -- if we include it, I think we need the equation plus (some variant of) the subsequent expository text}
%since
%\eqnn{
%$
%M \cdot \delta x_s^- = \delta x_s^- + \paren{F^+-F^-}\cdot \paren{\eta^\tr \cdot \delta x_s^-}/\paren{\eta^\tr \cdot F^-},
%$
%}
%
%\SB{variant 1:}
%which computation consists of two $d$-dimensional vector dot products, one $d$-dimensional vector scaling, and one $d$-dimensional vector addition.
%
%\SB{variant 2:}
%which computation only requires two $d$-dimensional vector dot products, avoiding the \Ord{d^2} time that would be required for multiplication by a general matrix.
%

% fig:pert
\begin{figure}[t]
\centering
\def\svgwidth{1.2\columnwidth} 
\resizebox{1.\columnwidth}{!}{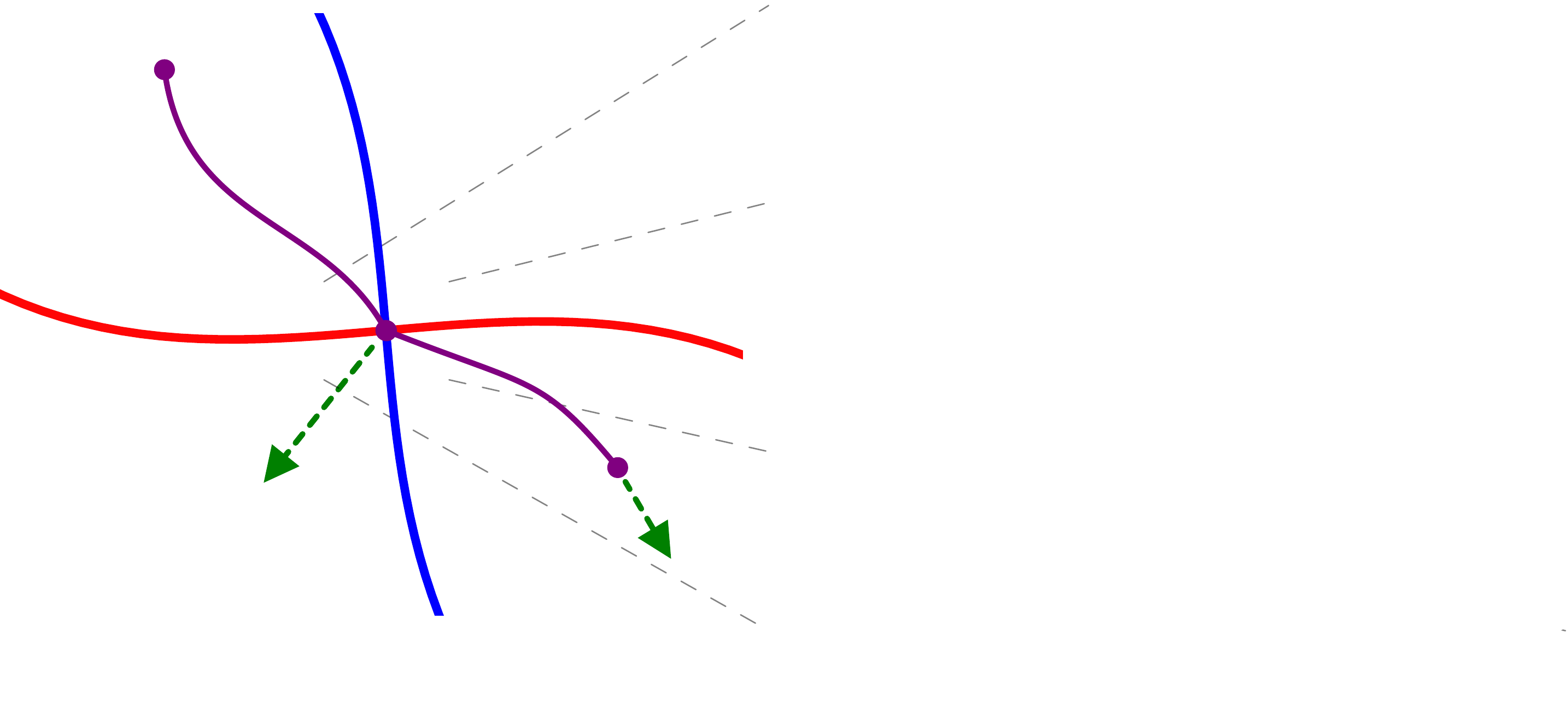}
\caption{%
\emph{Variational dynamics that determine the B-derivative of an {\ECr} vector field's {\PCr} flow~\eqref{eq:alg:B}.}
(a) Vector field $F:\R^2\into T\R^2$ is smooth everywhere except the smooth codimension-1 submanifolds $H_1, H_2 \subset\R^2$ that intersect transversally at $x_s\in\R^2$,
generating a piecewise-differentiable flow $\phi:[0,\infty)\times\R^2\into\R^2$ satisfying $\phi_\tau(x_0) = x_\tau$ for all $\tau\in[0,t]$, i.e.\ $F$ is {\ECr} and $\phi$ is {\PCr}~\cite{Burden2016-bb}.
The B-derivative $D \phi_t(x_0; \delta x_0) = \delta x_t$ is determined as in~\eqref{eq:alg:Dphi} by the continuous-time variational dynamics 
$\delta\dot{x}_\tau = DF(x_\tau)\cdot\delta x_\tau$ % in~\eqref{eq:alg:var}
and the discrete-time variational dynamics
$\delta x_s^+ = B(\delta x_s^-)$.
%where $\rho = x_s \in H_1 \cap H_2$. % determined by multiple applications of~\eqref{eq:alg:salt}.
The algorithms in~\figref{alg} evaluate the piecewise-linear function $B$ using the auxiliary nonsmooth system in (b) determined by the tangent planes $\samp{H}_1$, $\samp{H}_2$ and vector field limits $F_b(\rho)$ in~\eqref{eq:Gamma} for $b\in\set{(-\ones,(+1,-1),(-1,+1),+\ones} = \set{-1,+1}^2$.
}
\label{fig:pert}
\end{figure}

\subsection{Multiple surfaces of discontinuity}
\label{sec:alg:n}
If $F$ is continuously differentiable everywhere except a finite set of smooth \mbox{codimension-1} submanifolds $\set{H_j}_{j=1}^n$ that intersect the trajectory $x$ transversally at only one point $x_s$ (see~\figref{pert}(a) for an illustration when $n=2$),
$s\in(0,t)$,
we showed in~\cite[Eqn.~(65)]{Burden2016-bb} that
the discrete-time update~\eqref{eq:alg:salt} is applied once for each surface. 
However, the order in which the updates are applied, and the limiting values of the vector field used to determine each update's saltation matrix, depend on $\delta x_0$.
If the surfaces intersect transversally, there are $n!$ different saltation matrices determined by $2^n$ vector field values, so considering all update orders requires factorial time and space.
To make these observations precise and specify the notation employed in~\cref{fig:pert,fig:alg}, we formally define the class of nonsmooth vector fields considered in this paper~\cite[Defs.~1,~2]{Burden2016-bb}:
\begin{definition}{(event-selected $\Cr$ ($\ECr$) vector field)}\label{def:ECr}
A vector field $F:D\into T D$ 
defined on an open domain $D\subset\R^d$
is 
\emph{event-selected $\Cr$ with respect to $h\in\Cr(U,\R^n)$ at $\rho\in \R^d$} 
if $U\subset D$ is an open neighborhood of $\rho$ and:
\begin{enumerate}[leftmargin=.5cm]
  \item (event functions) there exists $f > 0$ such that $Dh(x) \cdot F(x) \ge f$ for all $x\in U$;
  \item (smooth extension) for all $b\in\set{-1,+1}^n = \Bn$, with
    \eqnn{
      D_b = \set{x\in U : b_j(h_j(x) - h_j(\rho)) \ge 0},
    }
    $F|_{\Int D_b}$ admits a $\Cr$ extension $F_b:U\into T U$.
\end{enumerate}
\end{definition}
Our algorithms in~\figref{alg}
compute
\eqnn{\label{eq:alg:B}
\delta x_s^+ = \delta\rho^+ = B(\delta\rho^-) = B(\delta x_s^-) %\in\R^d
}
given $\delta\rho^- = \delta x_s^-\in\R^d$, 
%compatibly-oriented%
%\footnote{
%We require $\eta^\tr\dot
%}
normals $\set{ \eta_j = Dh_j(\rho) }_{j=1}^n\subset\R^d$ at $x_s$
to surfaces $\set{ H_j = h_j^{-1}(\rho) }_{j=1}^n$%
%that are positively oriented%
%\footnote{
%We require that $\eta_j^\tr\cdot F > 0$ near $x_s$ as in~\cite[Def.~1,~2]{Burden2016-bb} to preclude
%\emph{sliding}~\cite{Utkin1977-rm, Jeffrey2014-nt} 
%and
%\emph{branching}~\cite[Def.~3.11]{Simic2005-fv}.
%}
%with respect to $F$%
,
and 
a function $\Gamma:\set{-1,+1}^n\into\R^d$ 
that evaluates limits of the vector field $F$ at $\rho = x_s$,
\eqnn{\label{eq:Gamma}
\forall b\in\set{-1,+1}^n : \Gamma(b) =  
F_b(\rho),
%\lim_{\alpha\goesto 0^+} F\paren{\rho + \alpha\sum_{j=1}^n b_j \eta_j},
}
using the piecewise-constant dynamics illustrated in~\figref{pert}(b), 
which are the discrete-time analog of the continuous-time variational dynamics~\eqref{eq:alg:var}.
Overall, the desired derivative is
\eqnn{\label{eq:alg:Dphi}
D \phi_t(x_0; \delta x_0) = D \phi_{t-s}(x_s) \cdot B\paren{D \phi_s(x_0) \cdot \delta x_0},
}
where $B:T_\rho\R^d\into T_\rho\R^d$ is the continuous piecewise-linear function defined by our algorithms in~\figref{alg}. 
Our algorithms require \Ord{n^2 d} time and \Ord{d} space to evaluate the directional derivative~\eqref{eq:dder}%
\footnote{%
These algorithms can be modified as in~\eqref{eq:comp:ind} to determine the order of surface crossings for the perturbed trajectory without changing the time or space complexity, so the associated saltation matrix~\eqref{eq:comp:Msigma} can be constructed in \Ord{n d^2} time and \Ord{d^2} space; this construction is discussed in more detail in~\secref{comp}.
}%
.

Assuming for the moment that these algorithms are correct, we emphasize that they achieve a dramatic reduction in the computational complexity of evaluating the B-derivative -- from factorial to low-order polynomial -- relative to {\naive} enumeration of all pieces of the B-derivative.
However, despite the apparent simplicity of our algorithms (computationally and conceptually), verifying their correctness requires significant technical effort; the bulk of the present paper is devoted to this verification task.

% alg:bdiff

\begin{figure}[t]
\begin{minipage}{0.46\textwidth}
\begin{algorithm}[H]
\caption{$\delta\rho^+ \assign B(\delta\rho^-,\eta,\Gamma)$}
\label{alg:B}
\begin{algorithmic}[1]
\STATE{$\delta t \assign 0\in\R$}
\STATE{$\delta\rho^+ \assign \delta\rho^-\in\R^d$}
\STATE{$b \assign -\ones\in\set{-1,+1}^n$}
\WHILE{$b \neq +\ones$}
\FOR{$j\in\set{1,\dots,n}$}%, b_j < 0$}
\STATE\label{alg:B:tti}{$\displaystyle \tau_j \assign - \paren{\eta^\tr_{j}\cdot\delta\rho^+} / \paren{\eta^\tr_{j}\cdot \Gamma(b)}$}
\ENDFOR
\STATE{$j^* \assign \arg\min_{j\in\set{1,\dots,n}}\set{ \tau_j : b_j < 0}$}
\STATE{$\delta t \assign \delta t + \tau_{j^*}$}
\STATE{$\delta \rho^+ \assign \delta \rho^+ + \tau_{j^*} \cdot \Gamma(b)$}
\STATE{$b_{j^*} \assign +1$}
\ENDWHILE
\RETURN $\delta\rho^+ - \delta t \cdot \Gamma(+\ones)$
\end{algorithmic}
\end{algorithm}
\vfill
\end{minipage}
\hfill
\renewcommand{\algorithmicdo}{\textbf{:}}
\begin{minipage}{0.51\textwidth}
\begin{algorithm}[H]
\caption{\tt def B(dx,e,G):}
\label{alg:PyB}
\begin{algorithmic}[1]
\ttfamily
%\STATE{import numpy as np}  % removed to make lines match
\STATE{dt = 0}
%\STATE{dx = np.asarray(dx).flatten()}
\STATE{dx = np.array(dx)}%.flatten()}
%\STATE{b = -np.ones(e.shape[0],dtype=np.int)}
\STATE{b = -np.ones(len(e),dtype=np.int)}
\WHILE{np.any(b < 0)}
\STATE\label{alg:PyB:tti}{tau = -np.dot(e,dx)/np.dot(e,G(b))}
\STATE{tau[b > 0] = np.inf}
\STATE{j = np.argmin(tau)}
\STATE{dt += tau[j]}
\STATE{dx += tau[j] * G(b)}
\STATE{b[j] = +1}
\ENDWHILE
\RETURN dx - dt * G(b) \# b ==  [+1,...,+1]
\end{algorithmic}
\end{algorithm}
\vfill
\end{minipage}
\caption{%
\emph{Algorithms that evaluate the B-derivative of an {\ECr} vector field's {\PCr} flow} 
written in pseudocode (\algref{B}) and {\tt Python}~\cite{Python_Software_Foundation_undated-lp} sourcecode (\algref{PyB}; requires {\tt import numpy as np}~\cite{Oliphant2006-ar}).
These algorithms apply at a point $\rho\in\R^d$ where a vector field $F:\R^d\into T\R^d$ is event-selected $C^r$ with respect to $n$ surfaces (see \figref{pert} for an illustration when $d=n=2$), and assume the following data is given:
}
\begin{minipage}{0.3\textwidth}
\footnotesize
\textit{
\begin{itemize}
  \item[ ] tangent direction,
  \item[ ] surface normals at $\rho$,
  \item[ ] vector field limits~\eqref{eq:Gamma},
\end{itemize}
}
\end{minipage}
\hfill
\begin{minipage}{0.30\textwidth}
\footnotesize
\textit{
\begin{itemize}
  \item[ ]$\delta\rho^-\in T_\rho \R^d$, 
  \item[ ]$\eta= \set{\eta_j}_{j=1}^n\subset\R^d$, 
  \item[ ]$\Gamma: \set{-1,+1}^n\into\R^d$, 
\end{itemize}
}
\end{minipage}
\hfill
\begin{minipage}{0.375\textwidth}
\footnotesize
\textit{
\begin{itemize}\normalfont
  \item[\texttt{dx}] -- \emph{array}, \texttt{dx.shape == (d,)};
  \item[\texttt{e}] -- \emph{array}, \texttt{e.shape == (n,d)};
  \item[\texttt{G}] -- \emph{function}, \texttt{G(b).shape == (d,)}.
\end{itemize}
}
\end{minipage}
\label{fig:alg}
\end{figure}

% Shai's edits
%\input{limbo}
%\input{bbdifflimbo}

\section{Applications}
\label{sec:app}

To illustrate and validate our methods, we apply the algorithm from the preceding section to 
%the class of
piecewise-constant vector fields 
%from~\cite[Sec.~8.1]{Burden2016-bb} 
in~\secref{app:pwc} 
and 
mechanical systems subject to unilateral constraints in~\secref{app:mech}.
Sourcecode implementation of \cref{alg:PyB} and applications from the remainder of this section are provided in SM.

%To assess the performance and accuracy of our proposed algorithm, we developed a Python 3.6.9 implementation utilizing the NumPy \cite{2020NumPy-Array} and SciPy \cite{2020SciPy-NMeth} numerical processing libraries.
%Our code available at:
%In order to validate both our algorithm and its implementation, we considered a suite of examples.
%
%If we posses explicit expressions for vector field $F : \mathbb{R}^d \to T \mathbb{R}^d$ and event surfaces $H_i, i \leq d$, the exact B-derivative components (and thus, the entire map) $M_{\sigma}$ can be computed exactly (i.e.\ symbolically) via \cite[Eqn. (67)]{Burden2016-bb}.
%For representatives  of this regime, we considered a first-order phase oscillator, \S 8.1 of \cite{Burden2016-bb}, a second-order phase oscillator of \S 8.2 of \cite{Burden2016-bb}, and perturbations of each.
%In the respective sections of \cite{Burden2016-bb}, symbolic expressions for $D_x \phi_{\omega}$ are derived, and we aim to reproduce the value of these expressions numerically via our algorithm given parameter values.

\subsection{Piecewise-constant vector field}
\label{sec:app:pwc}

Consider the vector field $F:\R^d\into T\R^d$ defined by
\eqnn{\label{eq:pwc}
  \dot{x} = F(x) = \ones + \Delta\paren{\sgn(x)}
}
%
%\begin{subequations}
%\label{eq:foso}
%\begin{align}
%  \dot{q} &= \alpha\cdot\ones - \delta\cdot\sgn(q) + \Delta_{\sgn(q)},
%  \label{eq:dq} \\
%  \ddot{p} &= \alpha\cdot\ones - \beta\cdot\dot{p} - \delta\cdot\sgn(p) + \Delta_{\sgn(q)}
%  \label{eq:ddp},
%\end{align}
%\end{subequations}
%
where 
$\Delta:\Bd\into\R^d$;
so long as all components of all vectors specified by $\Delta$ are larger than $-1$, i.e.\ $\min_{b\in\Bd}\brak{\Delta(b)}_j > -1$, $F$ is event-selected {\Cinfty} with respect to the identity function $h:\R^d\into\R^d$ defined by $h(q) = q$.
We regard~\eqref{eq:pwc} as a canonical form for piecewise-constant event-selected {\Cinfty} vector fields that are discontinuous across $d$ subspaces, since any such vector field can be obtained by applying a linear change-of-coordinates to~\eqref{eq:pwc}.
In what follows, we focus on the trajectory that passes through the origin $\rho = 0$, which lies at the intersection of $d$ surfaces of discontinuity for $F$.
With $\rho^- = \rho - \frac{1}{2} F_{-\ones}(\rho)$, $\rho^+ = \rho + \frac{1}{2} F_{+\ones}(\rho)$,
we note that $\rho^-$ flows to $\rho^+$ through $\rho$ in $1$ (one) unit of time.

\begin{figure}[t]
  \centering
  \resizebox{0.47\columnwidth}{!}{%% Creator: Inkscape 1.0.2 (e86c870879, 2021-01-15), www.inkscape.org
%% Accompanies image file '2d-sym.pdf' (pdf, eps, ps)
\begingroup%
  \makeatletter%
  \providecommand\color[2][]{%
    \errmessage{(Inkscape) Color is used for the text in Inkscape, but the package 'color.sty' is not loaded}%
    \renewcommand\color[2][]{}%
  }%
  \providecommand\transparent[1]{%
    \errmessage{(Inkscape) Transparency is used (non-zero) for the text in Inkscape, but the package 'transparent.sty' is not loaded}%
    \renewcommand\transparent[1]{}%
  }%
  \providecommand\rotatebox[2]{#2}%
  \newcommand*\fsize{\dimexpr\f@size pt\relax}%
  \newcommand*\lineheight[1]{\fontsize{\fsize}{#1\fsize}\selectfont}%
  \ifx\svgwidth\undefined%
    \setlength{\unitlength}{713.40362549bp}%
    \ifx\svgscale\undefined%
      \relax%
    \else%
      \setlength{\unitlength}{\unitlength * \real{\svgscale}}%
    \fi%
  \else%
    \setlength{\unitlength}{\svgwidth}%
  \fi%
  \global\let\svgwidth\undefined%
  \global\let\svgscale\undefined%
  \makeatother%
  \begin{picture}(1,0.90139373)%
    \lineheight{1}%
    \setlength\tabcolsep{0pt}%
    \put(0,0){\includegraphics[width=\unitlength,page=1]{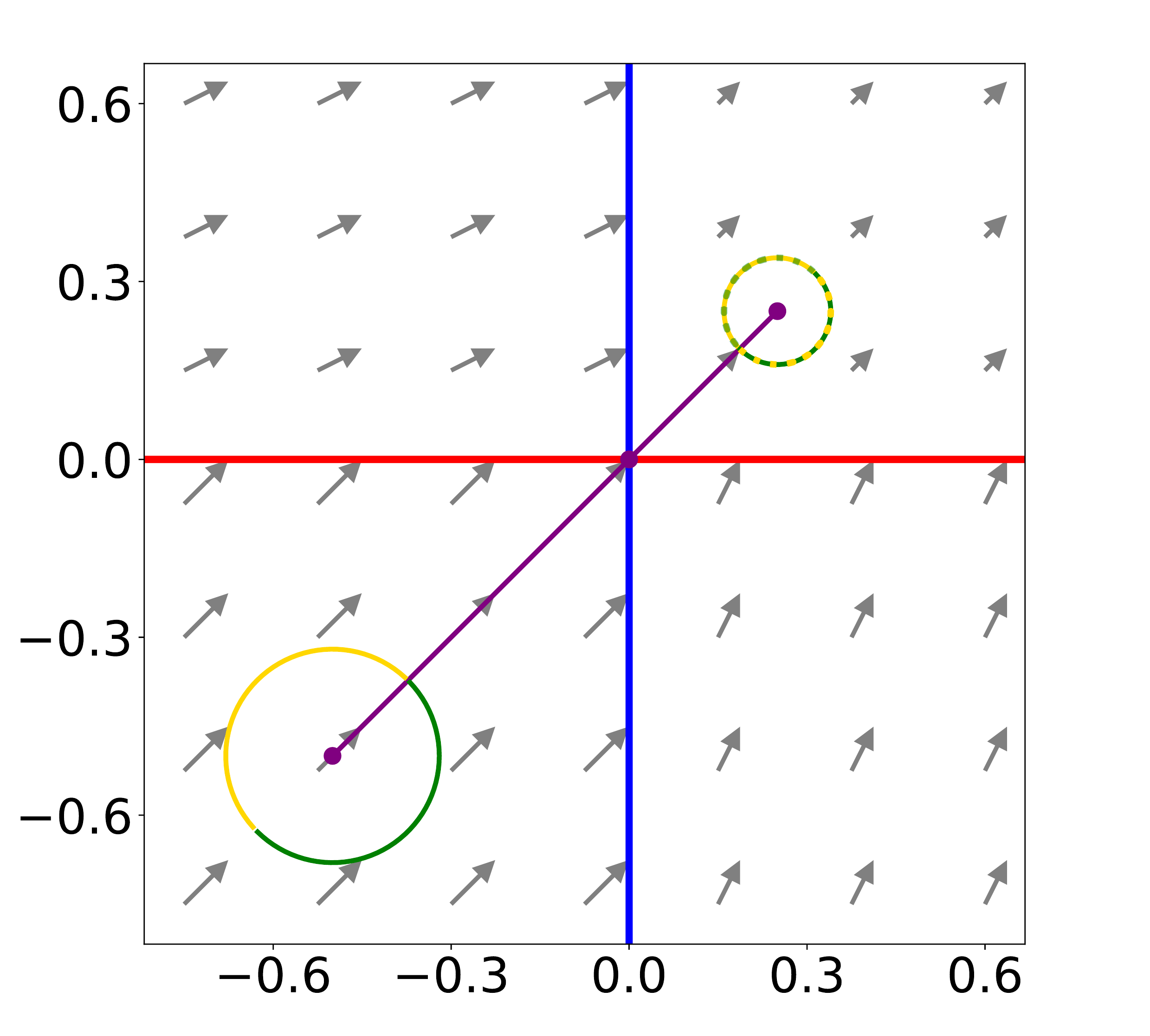}}%
    \put(0.13749788,0.75407364){\color[rgb]{0,0,0}\makebox(0,0)[lt]{\lineheight{1.25}\smash{\begin{tabular}[t]{l}\scalebox{3}{$F_{(-1,1)}=(1,\frac{1}{2})$}\end{tabular}}}}%
    \put(0.57228399,0.7594314){\color[rgb]{0,0,0}\makebox(0,0)[lt]{\lineheight{1.25}\smash{\begin{tabular}[t]{l}\scalebox{3}{$F_{(1,1)}=(\frac{1}{2},\frac{1}{2})$}\end{tabular}}}}%
    \put(0.13308053,0.41451361){\color[rgb]{0,0,0}\makebox(0,0)[lt]{\lineheight{1.25}\smash{\begin{tabular}[t]{l}\scalebox{3}{$F_{(-1,-1)}=(1,1)$}\end{tabular}}}}%
    \put(0.56512296,0.18536403){\color[rgb]{0,0,0}\makebox(0,0)[lt]{\lineheight{1.25}\smash{\begin{tabular}[t]{l}\scalebox{3}{$F_{(1,-1)}=(\frac{1}{2},1)$}\end{tabular}}}}%
    \put(0.55621623,0.42608093){\color[rgb]{0,0,0}\makebox(0,0)[lt]{\lineheight{1.25}\smash{\begin{tabular}[t]{l}\scalebox{3}{$H_1$}\end{tabular}}}}%
    \put(0.4627128,0.52073392){\color[rgb]{0,0,0}\makebox(0,0)[lt]{\lineheight{1.25}\smash{\begin{tabular}[t]{l}\scalebox{3}{$H_2$}\end{tabular}}}}%
  \end{picture}%
\endgroup%
}
  \resizebox{0.49\columnwidth}{!}{%% Creator: Inkscape 1.0.2 (e86c870879, 2021-01-15), www.inkscape.org
%% Accompanies image file '2d-asym.pdf' (pdf, eps, ps)
\begingroup%
  \makeatletter%
  \providecommand\color[2][]{%
    \errmessage{(Inkscape) Color is used for the text in Inkscape, but the package 'color.sty' is not loaded}%
    \renewcommand\color[2][]{}%
  }%
  \providecommand\transparent[1]{%
    \errmessage{(Inkscape) Transparency is used (non-zero) for the text in Inkscape, but the package 'transparent.sty' is not loaded}%
    \renewcommand\transparent[1]{}%
  }%
  \providecommand\rotatebox[2]{#2}%
  \newcommand*\fsize{\dimexpr\f@size pt\relax}%
  \newcommand*\lineheight[1]{\fontsize{\fsize}{#1\fsize}\selectfont}%
  \ifx\svgwidth\undefined%
    \setlength{\unitlength}{782.67504883bp}%
    \ifx\svgscale\undefined%
      \relax%
    \else%
      \setlength{\unitlength}{\unitlength * \real{\svgscale}}%
    \fi%
  \else%
    \setlength{\unitlength}{\svgwidth}%
  \fi%
  \global\let\svgwidth\undefined%
  \global\let\svgscale\undefined%
  \makeatother%
  \begin{picture}(1,0.83846161)%
    \lineheight{1}%
    \setlength\tabcolsep{0pt}%
    \put(0,0){\includegraphics[width=\unitlength,page=1]{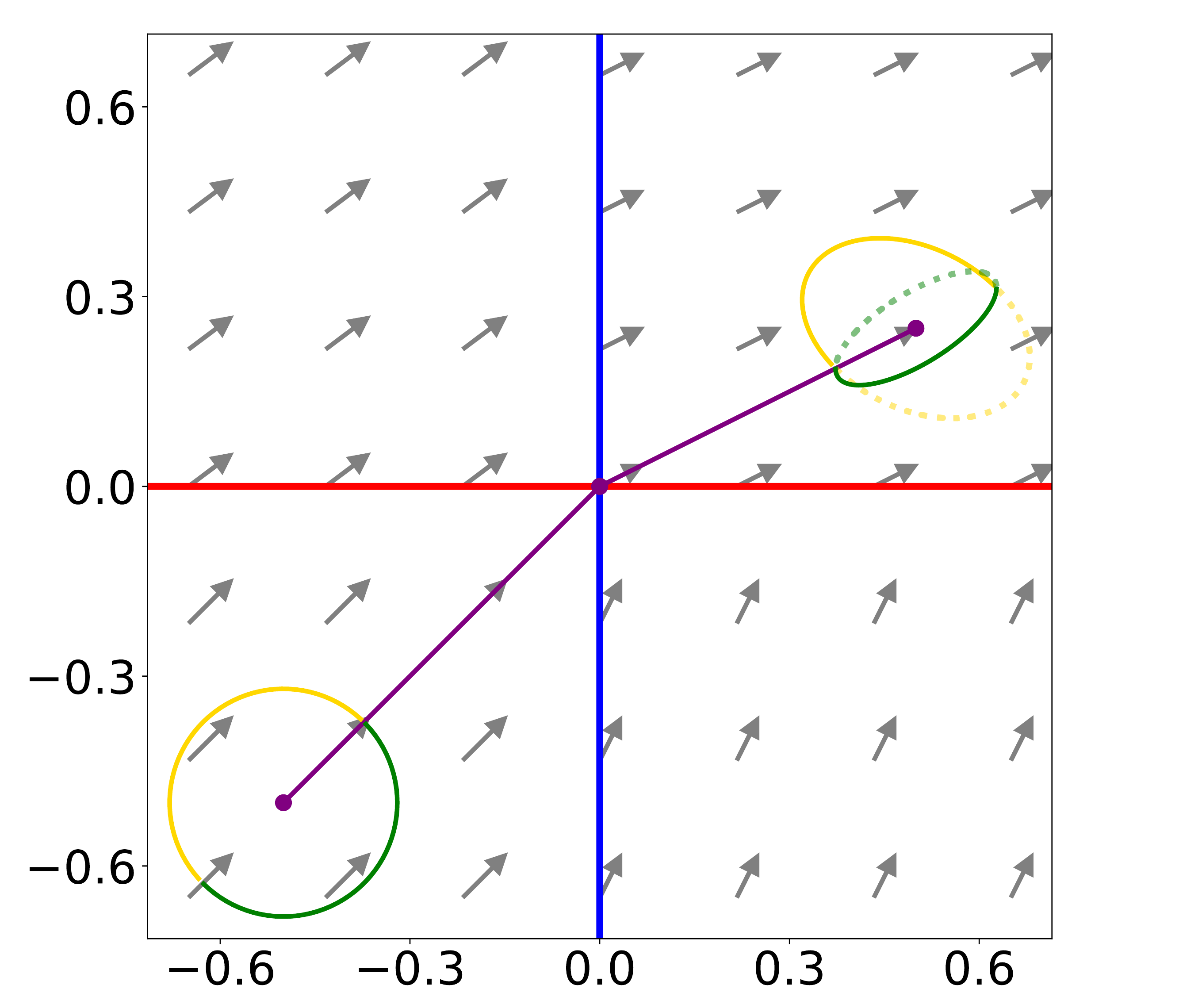}}%
    \put(0.13603821,0.72034666){\color[rgb]{0,0,0}\makebox(0,0)[lt]{\lineheight{1.25}\smash{\begin{tabular}[t]{l}\scalebox{3}{$F_{(-1,1)}=(1,\frac{3}{4})$}\end{tabular}}}}%
    \put(0.56389098,0.71979097){\color[rgb]{0,0,0}\makebox(0,0)[lt]{\lineheight{1.25}\smash{\begin{tabular}[t]{l}\scalebox{3}{$F_{(1,1)}=(1,\frac{1}{2})$}\end{tabular}}}}%
    \put(0.13527542,0.37385259){\color[rgb]{0,0,0}\makebox(0,0)[lt]{\lineheight{1.25}\smash{\begin{tabular}[t]{l}\scalebox{3}{$F_{(-1,-1)}=(1,1)$}\end{tabular}}}}%
    \put(0.5160252,0.15954482){\color[rgb]{0,0,0}\makebox(0,0)[lt]{\lineheight{1.25}\smash{\begin{tabular}[t]{l}\scalebox{3}{$F_{(1,-1)}=(\frac{1}{2},1)$}\end{tabular}}}}%
    \put(0.51008252,0.37786901){\color[rgb]{0,0,0}\makebox(0,0)[lt]{\lineheight{1.25}\smash{\begin{tabular}[t]{l}\scalebox{3}{$H_1$}\end{tabular}}}}%
    \put(0.20728333,0.44565112){\color[rgb]{0,0,0}\makebox(0,0)[lt]{\lineheight{1.25}\smash{\begin{tabular}[t]{l}\scalebox{3}{$H_2$}\end{tabular}}}}%
  \end{picture}%
\endgroup%
}
  \caption{\label{fig:fo}
  \emph{B-derivative of vector field from~\secref{app:pwc} in linear
  \emph{(left)} and piecewise-linear \emph{(right)} cases.}
  The vector field $F$ defined in~\eqref{eq:pwc} is piecewise-constant and discontinuous across the coordinate hyperplanes $H_1$, $H_2$,
  generating a piecewise-differentiable flow $\phi$ with B-derivative $B$.
  \emph{(left)} 
  The B-derivative is linear 
  in the special case defined by~\eqref{eq:pwc:lin}. 
  \emph{(right)} 
  The B-derivative is continuous and piecewise-linear in general, so a ball of initial conditions %(solid gold and green in lower-left quadrant) 
  flows to a piecewise-ellipsoid (gold and green fill). % in upper-right quadrant).
  }
\end{figure}

Our goal is to compute $D_x \phi(1,\rho^-;\delta\rho^-)\in T_{\rho^+}\R^d$ for a given $\delta\rho^-\in T_{\rho^-}\R^d$.
In the general case,
the desired derivative is piecewise-linear with (up to) $d!$ distinct pieces,
providing a general test.
In the special case where $\Delta(b) = -\delta\cdot b$ for all $b\in\Bd$, $\abs{\delta} < 1$, the desired derivative is linear~\cite[Eqn.~(86)]{Burden2016-bb},
\eqnn{\label{eq:pwc:lin}
%\paren{\set{\Delta_b}_{b\in\Bd} = \set{0_d}} \implies \paren
{D_x \phi(1,\rho^-;\delta\rho^-) = \frac{1 - \delta}{1 + \delta}\cdot\delta\rho^-},
}
providing a closed-form expression for comparison.
\figref{fo} illustrates results from both cases with $d = 2$;
a more exhaustive test suite is provided in SM. %~\cite{SM}.
%The supplementary materials contain a Python 3.6.9 implementation named ``test-bd.py'' which implements these two examples, along with related piecewise-constant cases. 

\subsection{Mechanical systems subject to unilateral constraints}
\label{sec:app:mech}

Consider a mechanical system whose configuration is subject to one-sided (i.e.\ \emph{unilateral}) constraints.
The dynamics of such systems have been studied extensively using 
the formalisms of
complementarity~\cite[Sec.~3]{Lotstedt1982-ea},
measure differential inclusions~\cite[Sec.~3]{Ballard2000-ui},
hybrid systems~\cite[Sec.~2.4,~2.5]{Johnson2016-nh},
and
geometric mechanics~\cite[Sec.~3]{Eldering2016-oj}.
Regardless of the chosen formalism, in a coordinate chart $Q\subset\R^d$ the dynamics governing $q$ take the form%
\footnote{%
We interpret the inequality $a(q) \geq 0$ componentwise. 
}
\eqnn{\label{eq:app:rigid}
M(q)\ddot{q} = f(q,\dot{q})\ \text{subject to}\ a(q) \geq 0
}
where:
$M(q)\in\R^{d\times d}$ specifies the kinetic energy metric;
$f(q,\dot{q})\in\R^d$ specifies the internal, applied, and Coriolis forces; %, where $u\in\R^m$ is a control input;
$a(q)\in\R^n$ specifies the \emph{unilateral constraints};
and
we assume in what follows that $M$, $f$, and $a$ are smooth functions.
Different formalisms enforce the constraint $a(q) \geq 0$ in~\eqref{eq:app:rigid} differently, so we consider several cases in the following subsections.
Additional illustrative examples are provided in SM.

\subsubsection{Rigid constraints yield discontinuous flows}
\label{sec:app:mech:rigid}
If constraints are enforced \emph{rigidly} as in~\cite{Lotstedt1982-ea, Ballard2000-ui, Johnson2016-nh}, 
meaning that they must be satisfied exactly,
then the velocity must undergo impact (i.e.\ change discontinuously) whenever $\dot{q}\in T_q Q$ is such that $a_j(q) = 0$ and $Da_j(q)\cdot \dot{q} < 0$ for some $j\in\set{1,\dots,n}$~\cite[Sec.~2]{Lotstedt1982-ea}~\cite[Eqn.~(23)]{Johnson2016-nh}~\cite[Eqn.~(23)]{Ballard2000-ui}.
Unfortunately for our purposes, these discontinuities in the state vector $x = (q,\dot{q})$ cannot be modeled using an event-selected {\Cr} vector field $\dot{x} = F(x)$,
and the flow of such systems is generally discontinuous%
\footnote{%
We note that the flow can be {\PCr} at non-impact times if the constraint surfaces intersect orthogonally~\cite{Pace2017-tt},
i.e. if the surface normals are orthogonal with respect to the inverse of the kinetic energy metric~\cite[Theorem~20]{Ballard2000-ui}. 
}%
.

\subsubsection{Soft conservative constraints yield Lipschitz vector fields, {\Co} flows}
\label{sec:app:mech:soft:Lip}
We now consider the formalism in~\cite{Eldering2016-oj} that ``softens'' (i.e.\ approximately enforces) rigid constraints $a(q) \geq 0$
by augmenting the potential energy with %continuously-differentiable 
\emph{penalty functions} $\set{v_j}_{j=1}^n$ that scale quadratically with the degree of constraint violation~\cite[Eqn.~(12)]{Eldering2016-oj},
\eqnn{
\forall j\in\set{1,\dots,n} : v_j(q) = \pw{ 0,& a_j(q) \geq 0 \\ \frac{1}{2}\kappa_j\, a_j^2(q),& a_j(q) < 0}
}
In essence, each rigid constraint $a_j(q) \geq 0$ is replaced by a spring with stiffness $\kappa_j$,
leading to the unconstrained dynamics~\cite[Eqn.~(14)]{Eldering2016-oj}
\eqnn{\label{eq:app:penalty}
M(q) \ddot{q} & = f(q,\dot{q},u) - \sum_{j=1}^n Dv_j(q)^\top \\
& = f(q,\dot{q},u) - \sum\set{\paren{\kappa_j\, a_j(q)}\cdot Da_j(q)^\top : j\in\set{1,\dots,n},\, a_j(q) < 0}.
}
As shown by~\cite[Thm.~3]{Tornambe1999-dr},
trajectories of~\eqref{eq:app:penalty} converge to those of~\eqref{eq:app:rigid} in the rigid limit (i.e.\ as stiffnesses go to infinity).
Importantly for our purposes, the dynamics in~\eqref{eq:app:penalty} can be modeled using an event-selected vector field %$\dot{x} = F(x)$ 
along trajectories that pass transversally through the constraint surfaces,
whence our algorithms can compute the B-derivative of the flow.
However, the vector field~\eqref{eq:app:penalty} in this case is (locally Lipschitz) continuous, hence the B-derivative is trivial (all non-identity terms in~\eqref{eq:comp:Msigma} are zero), whence the flow is continuously-differentiable ({\Co}).

\subsubsection{Soft dissipative constraints yield {\ECr} vector fields, {\Co} flows}
\label{sec:app:mech:soft:PCr}
We now augment the unconstrained dynamics~\eqref{eq:app:penalty} with dissipation as in~\cite{Eldering2016-oj}:
\eqnn{\label{eq:app:diss}
M(q) \ddot{q} 
& = f(q,\dot{q},u) - \sum\set{\paren{\kappa_j\, a_j(q) + \beta_j\, Da_j(q)\cdot\dot{q}}\cdot Da_j(q)^\top : j\in\set{1,\dots,n},\, a_j(q) < 0};
}
in essence, each constraint penalty is augmented by a spring-damper that is only active when the constraint is violated as in studies involving contact with complex geometry~\cite{Elandt2019-yz} or terrain~\cite{Aguilar2015-fy}.
The dynamics in~\eqref{eq:app:diss} can be modeled using an event-selected vector field along trajectories that pass transversally through the constraint surfaces, and the vector field is discontinuous along the constraint surfaces. %whence our algorithms can be applied to compute the (non-trivial) B-derivative.
However, we can show that the flow of~\eqref{eq:app:diss} is continuously-differentiable ({\Co}) along \emph{any} trajectory 
that passes transversally through constraint surfaces.
Indeed, letting $x = (q,\dot{q})$ denote the state of the system so that $\dot{x} = (\dot{q},\ddot{q}) = F(x)$ is determined by~\eqref{eq:app:diss},
the saltation matrix~\eqref{eq:alg:salt} associated with each constraint $a_j$ has the form
\eqnn{\label{eq:app:mech:salt}
%M_j = 
I + \frac{1}{Da_j(q)\cdot\dot{q}}\left[\begin{array}{c} 0 \\ \pm \paren{\kappa_j a_j(q) + \beta_j Da_j(q)\cdot\dot{q}} \cdot Da_j(q)^\top \end{array}\right]\left[\begin{array}{cc} Da_j(q) & 0 \end{array}\right]
}
where the sign in the column vector is 
determined by whether the constraint is activating ($-$) or deactivating ($+$).
Since matrices of the form in~\eqref{eq:app:mech:salt} commute, the saltation matrices associated with simultaneous activation and/or deactivation of multiple constraints are all equal, whence the flow of~\eqref{eq:app:diss} is continuously-differentiable ({\Co}) along any trajectory that passes transversally through constraint surfaces.

\begin{figure}[t]
  \centering
  % \includegraphics[width=0.5\textwidth]{chair-new.pdf}
  %\resizebox{\columnwidth}{!}{\input{figures/chair-new.pdf_tex}}
  %\includegraphics[width=\textwidth]{Figure_1.png}
  %\resizebox{\columnwidth}{!}{\input{figures/chair-big.pdf_tex}}
  \resizebox{\columnwidth}{!}{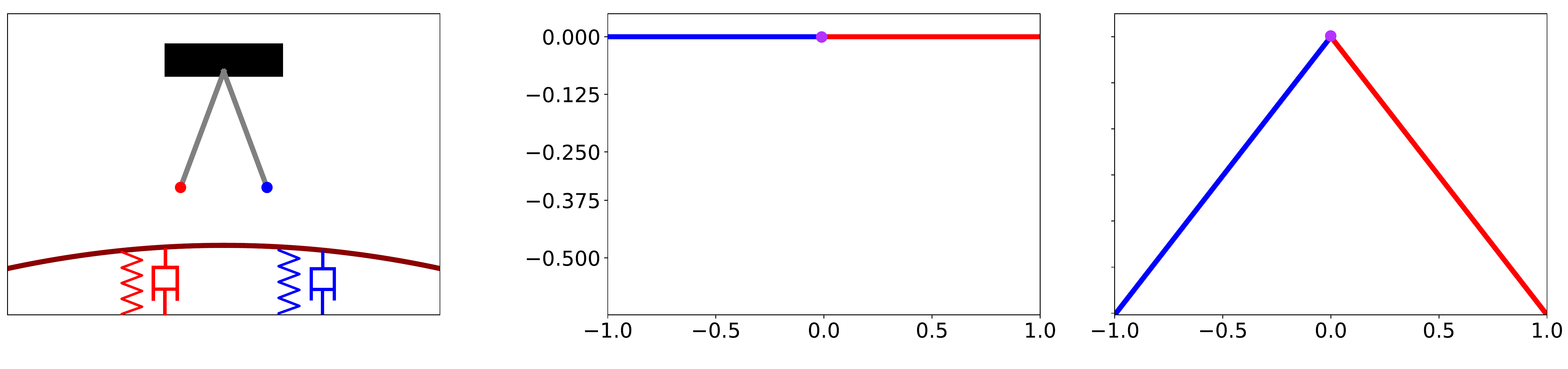}
  \caption{\label{fig:mech}
    \emph{Vertical-plane biped, a mechanical system subject to unilateral constraints (\secref{app:mech:biped}),}
    consists of a body with two rigid massless legs falling under the influence of gravity toward a substrate.
    The system's flow can be {\Co} (left column) or {\PCr} (right) depending on how forces vary as limbs contact substrate.
  }
\end{figure}

\subsubsection{Example (vertical-plane biped)}
\label{sec:app:mech:biped}
To ground the preceding observations, we consider the vertical-plane biped illustrated in~\figref{mech}(\emph{left}) that falls under the influence of gravity toward a substrate.
The biped body has mass $m$ and moment-of-inertia $J$;
we let $(x,y)\in\mathbb{R}^2$ denote the position of its center-of-mass in the plane and $\theta\in S^1$ denote its rotation. 
Two rigid massless limbs of length $\ell$ protrude at an angle of $\pm\psi$ with respect to vertical from the body's center-of-mass above a smooth substrate whose height is a quadratic function of horizontal position,
yielding unilateral constraints
\eqnn{\label{eq:app:mech:con}
a_1(x,y,\theta) = -y-\left(x+\ell \cos\left(\theta-\psi \right) \right)^2 - \ell \sin(\theta-\psi), \\
a_2(x,y,\theta) = -y-\left(x+\ell \cos\left(\theta+\psi \right) \right)^2 - \ell \sin(\theta+\psi).
}
We consider the smoothness of the system's flow along a trajectory that activates both constraints simultaneously%
\footnote{e.g.\ initial condition $\paren{(x_0,y_0,\theta_0),(\dot{x}_0,\dot{y}_0,\dot{\theta}_0)} = \paren{(0,h,0), (0,0,0)}$ where $h$ is the initial body height}%
.
Direct calculation%
\footnote{\label{fn:SM}%
Sourcecode that verifies this fact using a computer algebra system is provided in SM.} %~\cite{SM}.} 
shows that adopting the formalism in~\eqref{eq:app:diss} yields continuously-differentiable flow for this system as illustrated in~\figref{mech}(\emph{middle}).

To obtain a flow that is piecewise-differentiable but \emph{not} continuously-differentiable, 
%we introduce forcing in the biped from~\figref{mech} that depends non-additively on the set of active constraints.
%Specifically, 
we modify the damping coefficients in~\eqref{eq:app:diss} using the following logic%
\footnote{%
Although we introduce this logic purely for illustrative purposes, we note that non-trivial dependence of forcing on the set of active constraints could be implemented physically using clutches~\cite{Collins2015-sk} or actuators~\cite{Sreenath2011-lb}.%
}%
:
$\beta_1 = \beta_2 = \frac{1}{2} $ if $a_1(q) < 0$ or $a_2(q) \geq 0$ (\emph{exclusive} or);  
%$\beta_2 = \frac{1}{2} $ if $a_2(q) < 0$ and $a_1(q) \geq 0$; 
$\beta_1 = \beta_2 = 1$ if $a_1(q) < 0$ and $a_2(q) <0$.
Direct calculation$^{\ref{fn:SM}}$ shows that the saltation matrices obtained from different sequences of constraint activations (left foot reaches substrate before right foot or vice-versa) are distinct:
\eqnn{\label{eq:app:mech:biped:salt}
M_{(\text{\color{red}{left}},\text{\color{blue}{right}})}
-
M_{(\text{\color{red}{right}},\text{\color{blue}{left}})}
=
\begin{bmatrix}  
	 0 &0&0 &0&0&0 \\
	 0 &0&0 &0&0&0 \\
	 0 &0&0 &0&0&0 \\
	 0 &0&0 &0&0&0 \\
	 -4 \beta \cos(\psi) & 0 & -2 \beta(\sin(2\psi)+\cos(\psi)) & 0 &0 &0 \\
	 0 &0&0 &0&0&0 \\
\end{bmatrix}.
}
The piecewise-linear B-derivative of the system's flow is illustrated in~\figref{mech}(\emph{right}).

\section{Background}
\label{sec:bg}

To verify correctness of the algorithms specified in~\secref{alg}, 
we utilize
the representation of piecewise-affine functions from~\cite{Groff2003-pb},
elements of the theory of piecewise-differentiable functions from~\cite{Scholtes2012-la},
and
results about the class of nonsmooth flows under consideration from~\cite{Burden2016-bb}.
In an effort to make this paper self-contained (i.e.\ to save the reader from needing to cross-reference multiple citations to follow our derivations), we include a substantial amount of background details in this section.
The expert reader may wish to skim or skip this section, returning only if questions arise in subsequent sections.

\subsection{Polyhedral theory}
\label{sec:bg:poly}
We let 
$0_d\in\Real^d$ denote the vector of zeros,
$\ones_n\in\R^n$ the vector of ones,
and
$I_d\in\Real^{d\times d}$ the identity matrix;
when dimensions are clear from context, we suppress subscripts.
The vectorized signum function $\sgn:\R^d\into\set{-1,+1}^d$ is defined by
\eqnn{
\forall x\in\R^d,\ j\in\set{1,\dots,d} : \brak{\sgn(x)}_j = \sgn(x_j) = \pw{-1, & x_j < 0;\\ +1,& x_j \geq 0}.
}
%\GC{I think the square brackets look weird, but I don't have an immediate improvement.}
%\SB{in~\cite{Burden2016-bb} we did $e_j^\tr\cdot \sgn(x) = \sgn(x_j)$, but I don't love that -- especially because we don't use $e_j$ anywhere else in this paper}
If $A\in\R^{\ell\times m}$ and $B\in\R^{m\times n}$ then $A\cdot B\in\R^{\ell\times n}$ denotes matrix multiplication.
Given a subset $S\subset\R^d$, we define~\cite[Sec.~2.1.1]{Scholtes2012-la} 
%\GC{this is a lot of overhead -- we should make super sure each one is used, and if not, deleted here}
%\SB{agreed -- currently $\aff$, $\cone$, $\conv$ are used and $\lin$ is not}

\begin{subequations}
\begin{align}
%\lin S & = \set{\sum_{j=1}^n \alpha_j\, v_j : n\in\N,\ \set{v_j}_{j=1}^n \subset S,\ \set{\alpha_j}_{j=1}^n \subset\R},\\
\aff S & = \set{\sum_{j=1}^n \alpha_j\, v_j : n\in\N,\ \set{v_j}_{j=1}^n \subset S,\ \set{\alpha_j}_{j=1}^n \subset\R,\ \sum_{j=1}^n \alpha_j = 1},\\
\cone S & = \set{\sum_{j=1}^n \alpha_j\, v_j : n\in\N,\ \set{v_j}_{j=1}^n \subset S,\ \set{\alpha_j}_{j=1}^n \subset[0,\infty)},\\
\conv S & = \set{\sum_{j=1}^n \alpha_j\, v_j : n\in\N,\ \set{v_j}_{j=1}^n \subset S,\ \set{\alpha_j}_{j=1}^n \subset[0,1],\ \sum_{j=1}^n \alpha_j = 1},
\end{align}
\end{subequations}
termed the 
%\emph{linear span}, % <-- Scholtes term: ``linear hull''
\emph{affine span}, % <-- Scholtes term: ``affine hull''
\emph{cone span}, % <-- Scholtes term: ``cone generated by''
and 
\emph{convex hull} 
of $S$, respectively.
The \emph{dimension} of a convex set $S$ is defined to be
the dimension of its affine span, $\dim S = \dim \aff S$.
% polyhedron
A nonempty set $S\subset\R^d$ is called a \emph{polyhedron}~\cite[Sec.~2.1.2]{Scholtes2012-la} if there exists $A\in\R^{m\times d}$, $b\in\R^m$ such that
$S = \set{x\in\R^d : A \cdot x \leq b}$;
note that $S$ is closed and convex.
The linear subspace $\Li = \set{x\in\R^d : A \cdot x = 0}$ is called the \emph{lineality space} of $S$.
%Given a finite collection of polyhedra $\Sigma = \set{\Sigma_\omega}_{\omega\in\Omega}$ in $\R^d$,
%the union $\car{\Sigma} = \bigcup_{\omega\in\Omega} \Sigma_\omega\subset\R^d$ is called the \emph{carrier} of $\Sigma$.
%
% [ ] carrier
% [ ] polyhedral subdivision 
% [ ] localization
% TODO cover these in subsequent section !!!
%
%we will call $\Sigma$ a \emph{polyhedral partition}%
%\footnote{%
%We deliberately avoid the term \emph{polyhedral subdivision}~\cite[Sec.~2.2.1]{Scholtes2012-la} since the collections of polyhedra we consider below do not necessarily partition a \emph{convex} whole; see~\figref{tri} for an illustration. % of convex and non-convex polyhedral partitions.
%}
%if (i) all polyhedra in $\Sigma$ have the same dimension and (ii) the intersection of any two polyhedra in $\Sigma$ is either empty or a proper face of both polyhedra.

\subsection{Piecewise-affine functions}
\label{sec:bg:pwa}

%We will say a function $P:D\into C$ where $D\subset\R^d$ and $C\subset\R^c$ is \emph{piecewise-affine} if (i) $P$ is continuous and (ii) there are a finite collection $\Sigma$ of subsets of $\R^d$ such that $\bigcup_{S\in\Sigma} S = D$ and $P|_S$ is affine for each $S\in\Sigma$.
%There are multiple ways to represent such functions; for instance, the definition suggests the representation
%\eqnn{\label{eq:bg:pwa}
%\forall S\in\Sigma, x\in S : P|_S(x) = A_S\cdot x + b_S.
%}
%However, the representation in~\eqref{eq:bg:pwa} suppresses the constraints imposed on (ii) \GC{grammar seems weird} the ``pieces'' $\set{P|_S}_{S\in\Sigma}$ by (i) the continuity of $P$ -- namely, that $P|_{S}(x) = P|_{S'}(x)$ for all $x\in S\cap S'$.
%\GC{Not sure this is enough if the sets S are allowed to be arbitrary. E.g, what if S = Q, S' = irrationals. Let P(S) = 1, P'(S')=0. Both seem fine for piecewise affine, but this function is the dirichlet function and is continuous nowhere. Need something like e.g., the S are open(or closed), and use the gluing lemma from point set topology}

% triangulation
We will represent a piecewise-affine function
using a \emph{triangulation} 
$(Z^-, Z^+, \Delta)$~\cite[Sec.~3.1]{Groff2003-pb}
that consists of 
a combinatorial simplicial complex $\Delta$ whose vertex set is in 1-to-1 correspondence with each of the finite sets of vectors $Z^-\subset\R^d$, $Z^+\subset\R^c$. 
For our purposes,%
\footnote{
There are more general definitions of ([complete] semi-)simplicial complexes and the closely-related concept of $\Delta$-complexes in the literature~\cite[Ch.~2.1]{Hatcher2002-ec},~\cite[App.~A.3.1]{Groff2003-pb}.
Since we employ these concepts primarily in service of parameterizing piecewise-affine functions as in~\cite[Sec.~3.1]{Groff2003-pb},
we adopt the (relatively restrictive) definitions of \emph{combinatorial} and \emph{geometric} simplicial complexes from~\cite[Sec.~2.2.1]{Groff2003-pb} in what follows.
}
% combinatorial simplicial complex
a \emph{combinatorial simplicial complex} $\Delta$ is a collection of finite sets $\Delta = \set{\Delta_\omega}_{\omega\in\Omega}$ 
such that
$S\subset\Delta_\omega \implies S\in\Delta$ for all $\omega\in\Omega$;
we call $\bigcup_{\omega\in\Omega} \Delta_\omega$ the \emph{vertex set} of $\Delta$.
% geometric simplicies
We assume that, for every $\omega\in\Omega$, the collections of vectors $Z_\omega^\pm\subset Z^\pm$ determined by $\Delta_\omega$ are \emph{affinely independent}~\cite[Sec.~2.1.1]{Groff2003-pb} 
so that $\Delta_\omega^\pm = \conv Z_\omega^\pm$ are $\paren{\card{\Delta_\omega}-1}$-dimensional geometric simplices~\cite[Claim~2.9]{Groff2003-pb} where $\Delta_\omega^-\subset\R^d$, $\Delta_\omega^+\subset\R^c$.
% geometric simplicial complex
We assume further that, for every $\omega,\omega'\in\Omega$, the collections of vectors $Z_{\omega,\omega'}^\pm\subset Z^\pm$ determined by $\Delta_\omega\cap\Delta_{\omega'}$ coincide with $Z_\omega^\pm\cap Z_{\omega'}^\pm\subset Z^\pm$
so that
$\Delta^\pm = \set{\Delta_\omega^\pm}_{\omega\in\Omega}$ are 
\emph{geometric simplicial complexes}~\cite[Sec.~2.2.1]{Groff2003-pb}.
% TODO observe that simplices are polyhedra?
%Each geometric simplex $\Delta_\omega^\pm$ is, in particular, a (compact) \emph{polyhedron} as in~\secref{bg:poly}. 
With these assumptions in place,
the correspondence between $Z^-$ and $Z^+$ determined by the triangulation $(Z^-,Z^+,\Delta)$ 
%defines $f|_{Z^-}$
%With these assumptions in place, the triangulation $(Z^-,Z^+,\Delta)$ 
uniquely defines a piecewise-affine function $P:\car{\Delta^-}\into\car{\Delta^+}$
using the construction from~\cite[Sec.~3.1]{Groff2003-pb}
where
$\car{\Delta^-} = \bigcup_{\omega\in\Omega}\Delta_\omega^-\subset\R^d$, $\car{\Delta^+} = \bigcup_{\omega\in\Omega}\Delta_\omega^+\subset\R^c$
are termed the \emph{carriers}~\cite[Sec.~2.2.1]{Scholtes2012-la} 
of the geometric simplicial complexes $\Delta^\pm$.

\subsection{Piecewise-linear functions}
\label{sec:bg:PWL}

If a piecewise-affine function $P:\R^d\into\R^c$ is \emph{positively homogeneous}, that is,
\eqnn{\label{eq:bg:hom}
\forall \alpha \geq 0, v\in\R^d : P(\alpha\cdot v) = \alpha\cdot P(v),
}
then $P$ is \emph{piecewise-linear}~\cite[Prop.~2.2.1]{Scholtes2012-la}.
In this case, $P$ admits a \emph{conical subdivision}~\cite[Prop.~2.2.3]{Scholtes2012-la},
that is,
there exists a finite collection
$\Sigma = \set{\Sigma_\omega}_{\omega\in\Omega}$ 
such that:
(i) $\Sigma_\omega\subset\R^d$ is a $d$-dimensional \emph{polyhedral cone} for each $\omega\in\Omega$;%
\footnote{i.e.\ $\Sigma_\omega = \set{\sum_{j=1}^{\ell_\omega} \alpha_j v^\omega_j : \set{\alpha_j}_{j=1}^{\ell_\omega}\subset [0,\infty)}$, some $\set{v_j}_{j=1}^{\ell_\omega}\subset\R^d$~\cite[Thm.~2.1.1]{Scholtes2012-la}, 
and $\dim\Sigma_\omega = d$}
(ii) the $\Sigma_\omega$'s cover $\R^d$;% 
\footnote{i.e.\ $\bigcup_{\omega\in\Omega} \Sigma_\omega = \R^d$}
and
(iii) the intersection $\Sigma_\omega\cap\Sigma_{\omega'}$ is either empty or a \emph{proper face} of both polyhedral cones for each $\omega,\omega'\in\Omega$.%
\footnote{i.e.\ $\Sigma_\omega\cap\Sigma_{\omega'} = \set{\sum_{j=1}^{\ell_{\omega,\omega'}} \alpha_j v^{\omega,\omega'}_j : \set{\alpha_j}_{j=1}^{\ell_{\omega,\omega'}}\subset [0,\infty)}$, 
some
$\set{v^{\omega,\omega'}_j}_{j=1}^{\ell_{\omega,\omega'}}\subset\set{v^\omega_j}_{j=1}^{\ell_\omega}\cup\set{v^{\omega'}_j}_{j=1}^{\ell_{\omega'}}$}

\subsection{Piecewise-differentiable ($\PCr$) functions}
\label{sec:bg:PCr}

(This section is largely repeated from~\cite[Sec.~3.2]{Burden2016-bb}.)
The notion of piecewise--differentiability we employ was originally introduced in~\cite{Robinson1987-va}; since the monograph~\cite{Scholtes2012-la} provides a more recent and comprehensive exposition, we adopt the notational conventions therein.
Let $r \in \N \cup \set{\infty}$ and $D \subset \R^d$ be open.
A continuous function $f : D \into \R^c$ is called \emph{piecewise-\Cr} if for every $x_0 \in D$ there exists an open set $U \subset D$ containing $x_0$ and a finite collection $\set{f_j : U \into \R^c }_{j\in\J}$ of $\Cr$ functions such that for all $x \in U$ we have $f(x) \in \set{f_j(x)}_{j\in\J}$.
The functions $\set{f_j}_{j\in\J}$ are called \emph{selection functions} for $f|_U$, and $f$ is said to be a \emph{continuous selection} of $\set{f_j}_{j\in\J}$ on $U$.
A selection function $f_j$ is said to be \emph{active} at $x\in U$ if $f(x) = f_j(x)$.
We let $\PCr(D,\R^c)$ denote the set of piecewise-$\Cr$ functions from $D$ to $\R^c$.
Note that {\PCr} is closed under composition.
%Any $f\in\PCr(D,\R^c)$ is locally Lipschitz continuous, and a Lipschitz constant for $f$ is given by the supremum of the induced norms of the (\Frechet) derivatives of a set of selection functions for $f$.
The definition of piecewise-\Cr may at first appear unrelated to the intuition that a function ought to be piecewise-differentiable precisely if its ``domain can be partitioned locally into a finite number of regions relative to which smoothness holds''~\cite[Section~1]{Rockafellar2003-nd}.
However, as shown in~\cite[Thm.~2]{Rockafellar2003-nd}, piecewise-\Cr functions are always piecewise-differentiable in this intuitive sense.

Piecewise-differentiable functions possess a first--order approximation $Df:TD\into T\R^c$ called the \emph{Bouligand derivative} (or B--derivative)~\cite[Ch.~3]{Scholtes2012-la}; this is the content of~\cite[Lemma~4.1.3]{Scholtes2012-la}.
Significantly, %and unlike the directional derivative, 
this B--derivative obeys generalizations of many techniques familiar from calculus, including
the Chain Rule~\cite[Thm~3.1.1]{Scholtes2012-la}, % (and hence Product and Quotient Rules~\cite[Cor.~3.1.1]{Scholtes2012-la}),
Fundamental Theorem of Calculus~\cite[Prop.~3.1.1]{Scholtes2012-la},
and Implicit Function Theorem~\cite[Cor.~20]{Ralph1997-yr}. 
We let $Df(x;\delta x)$ denote the B--derivative of $f$ evaluated on the tangent vector $\delta x \in T_x D$.
The B-derivative is positively homogeneous, i.e. $\forall \delta x\in T_x D,\lambda\ge 0 : Df(x;\lambda\, \delta x) = \lambda Df(x; \delta x)$, and coincides with the directional derivative of $f$ in the $\delta x\in T_x D$ direction.
In addition, the B-derivative $Df(x):T_x D\into T_{f(x)}\R^c$ of $f$ at $x\in D$ is a continuous selection of the derivatives of the selection functions active at $x$~\cite[Prop.~4.1.3]{Scholtes2012-la},
\eqnn{
\forall \delta x \in T_x D : Df(x;\delta x)\in\set{Df_j(x)\cdot \delta x}_{j\in \J}.
}
However, the function $Df$ is generally \emph{not} continuous at $(x,\delta x)\in T D$; if it is, then $f$ is $C^1$ at $x$~\cite[Prop.~3.1.2]{Scholtes2012-la}.

\subsection{Event-selected {\Cr} ({\ECr}) vector fields and their {\PCr} flows}
\label{sec:bg:ECr}

Vector fields with discontinuous right-hand-sides and their associated flows have been studied extensively~\cite{Filippov1988-nh}.
In~\cref{def:ECr}~\cite[Defs.~1,~2]{Burden2016-bb}, a special class of so-called \emph{event-selected \Cr} (\ECr) vector fields 
were defined which
are allowed to be discontinuous along a finite number of codimension-1 submanifolds 
%-- the level sets $H_j = h_j^{-1}(h_j(\rho))$ --
but do not exhibit \emph{sliding}~\cite{Jeffrey2014-nt} along these submanifolds, and are $\Cr$ elsewhere.
Importantly, as shown in~\cite[Thm.~5]{Burden2016-bb}, an event-selected $\Cr$ vector field $F:\R^d\into T\R^d$ generates a piecewise-differentiable 
flow, that is, there exists a function $\phi:\F\into \R^d$ 
that is \emph{piecewise-$\Cr$} ($\phi\in\PCr$) in the sense defined in~\cite[Sec.~4.1]{Scholtes2012-la} (summarized in~\secref{bg:PCr})
where $\F\subset\R\times\R^d$ and 
\eqnn{
\forall (t,x)\in\F : \phi(t,x) = x + \int_0^t F(\phi(s,x)) ds.
}
Since $\phi$ is $\PCr$, it admits a first-order approximation $D\phi:T\F\into T\R^d$ termed the \emph{Bouligand} (or \emph{B-})derivative~\cite[Sec.~3.1]{Scholtes2012-la}, which is a continuous piecewise-linear function of tangent vectors at every $(t,x)\in\F$, that is,
the directional derivative $D\phi(t,x):T_{(t,x)}\F\into T_{\phi(t,x)}\R^d$ is continuous and piecewise-linear for all $(t,x)\in\F$.

\subsection{B-derivative of an {\ECr} vector field's {\PCr} flow}
\label{sec:bg:B}

Suppose $F:\R^d\into T\R^d$ is an {\ECr} vector field with {\PCr} flow $\phi:\F\into\R^d$.
Given a tangent vector $(\delta t,\delta x)\in  T_{(t,x)}\F$, 
it was shown in~\cite[Sec.~7.1.4]{Burden2016-bb} that
the value of the B-derivative $ D\phi(t,x;\delta t,\delta x)\in T_{\phi(t,x)} \R^d$ can be obtained by solving a jump-linear-time-varying differential equation~\cite[Eqn.~(70)]{Burden2016-bb},
where the ``jump'' arises from a matrix $\Xi_\omega$ determined by the sequence $\omega$ in which the perturbed initial state $x + \alpha\,\delta x$ crosses the surfaces of discontinuity of the vector field $F$ for small $\alpha > 0$~\cite[Eqn.~(67)]{Burden2016-bb}.
However,~\cite{Burden2016-bb} did not provide a representation of the piecewise-linear operator $D\phi(t,x)$ (and, to the best of our knowledge, neither has subsequent work).
The key theoretical contribution of this paper, obtained in~\secref{rep}, is a representation of the B-derivative with respect to state, $D_x\phi(t,x)$, using a \emph{triangulation} of its domain and codomain as defined in~\cite[Sec.~3.1]{Groff2003-pb} (and recalled in~\secref{bg:pwa}).

To inform the triangulation of the B-derivative $D_x\phi(t,x)$, we recall the values it takes on.
Since the flow $\phi:\F\into \R^d$ is piecewise-$\Cr$ ($\PCr$), it is a \emph{continuous selection} of a finite collection of $\Cr$ functions $\set{\phi_\omega:\F_\omega\into \R^d}_{\omega\in\Omega}$ near $(t,x)\in\F$, where $\F_\omega\subset\F$ is an open set containing $(t,x)$ for each $\omega\in\Omega$~\cite[Sec.~4.1]{Scholtes2012-la},
and
the B-derivative $D_x \phi(t,x)$ is a continuous selection of the classical (\emph{\Frechet} or \emph{F-})derivatives $\set{D_x\phi_\omega(t,x)}_{\omega\in\Omega}$~\cite[Prop.~4.1.3]{Scholtes2012-la},
that is,
\eqnn{\label{eq:Dxphi}
  \forall\delta x \in W_\omega\subset T_x \R^d : 
  D_x\phi(t,x;\delta x)
  & = 
  D_x \phi_\omega(t,x)\cdot\delta x,
}
where $W_\omega\subset T_x \R^d$ is the subset of tangent vectors where the selection function $D_x \phi_\omega$ is \emph{essentially active}~\cite[Prop.~4.1.1]{Scholtes2012-la}.
If 
$s,t\in\R$ 
and 
$x\in \R^d$ 
are such that 
$0 < s < t$ 
and 
the vector field $F$ is $\Cr$ on $\phi([0,t]\sm\set{s},x)$,
i.e. the trajectory initialized at $x\in \R^d$ encounters exactly one discontinuity of $F$ at
$\rho = \phi(s,x)$
on the time interval $[0,t]$,
then $D_x\phi_\omega(t,x)$ has the form
\eqnn{\label{eq:Dxphiomega}
  D_x \phi_\omega(t,x)
  = 
  D_x \phi(t-s,\rho) 
  \cdot
  \begin{bmatrix} F_{+\ones}(\rho) & I_d \end{bmatrix}
  \cdot
  \Xi_{\omega}
  \cdot
  \begin{bmatrix} 0_d^\top \\ I_d \end{bmatrix} 
  \cdot
  D_x \phi(s,x)
}
where
$F_{+\ones}$ is the $\Cr$ extension of $F|_{\Int{D_{+\ones}}}$ that exists by virtue of condition 2 in Def.~\ref{def:ECr}
and
$\Xi_\omega \in \Real^{(d+1)\times(d+1)}$ is the matrix from~\cite[Eqn.~(67)]{Burden2016-bb} corresponding to the selection function index $\omega\in\Omega$.
In what follows, we will work in circumstances where the selection functions are indexed by the symmetric permutation group over $n$ elements, i.e. $\Omega = \Sn$, and combine~\eqref{eq:Dxphi} and~\eqref{eq:Dxphiomega} as
\eqnn{\label{eq:Dxphisigma}
  \forall\delta x\in W_\sigma\subset T_x \R^d 
  :
  D_x \phi(t,x; \delta x)
  = 
  D_x \phi(t-s,\rho) 
  \cdot
  M_\sigma
  \cdot
  D_x \phi(s,x)\cdot\delta x
}
where the \emph{saltation matrix}%
\footnote{$\Xi_\sigma\in\Real^{(d+1)\times(d+1)}$ is referred to as a saltation matrix in~\cite[Sec.~7.1.4]{Burden2016-bb}, but this usage is inconsistent with the original definition of $M_\sigma\in\Real^{d\times d}$ as the saltation matrix in~\cite{Aizerman1958-ih}.}
$M_\sigma\in\Real^{d\times d}$ corresponding to index $\sigma$ is defined by
\eqnn{\label{eq:Msigma}
  M_\sigma 
  = 
  \begin{bmatrix} F_{+\ones}(\rho) & I_d \end{bmatrix} 
  \cdot
  \Xi_{\sigma}
  \cdot
  \begin{bmatrix} 0_d^\top \\ I_d \end{bmatrix}.
}

\subsection{Local approximation of an {\ECr} vector field}
\label{sec:bg:samp}

Suppose vector field $F:\R^d\into T\R^d$ is event-selected $\Cr$ with respect to $h\in\Cr(U,\R^n)$ at $\rho\in U\subset\R^d$.
For $b\in\Bn = \set{-1,+1}^n$ let
\eqnn{
\samp{D}_b = \set{x\in\R^d : b_j\, Dh_j(\rho)(x - \rho) \ge 0}
}
and consider the piecewise-constant vector field $\samp{F}:\R^d\into T\R^d$ defined by
\eqnn{
\forall b\in\Bn,\ x\in\samp{D}_b : \samp{F}(x) = F_b(\rho)
}
where $F_b$ is the $\Cr$ extension of $F|_{\Int D_b}$ that exists by virtue of condition~2 in Def.~\ref{def:ECr}\footnote{Note that $\samp{F}$ is well-defined as the value of $F_b$ is uniquely determined at $\rho$ by virtue of being continuous, even though the original $F$ is undefined at $\rho$.}
Note that $\samp{F}$ is event-selected $\Cr$ with respect to the affine function $\samp{h}$ defined by
\eqnn{
\forall x\in\R^d : \samp{h}(x) = Dh(\rho)(x-\rho),
}
whence it generates a piecewise-differentiable flow $\samp{\phi}:\samp{\F}\into\R^d$ where $\samp{\F} = \R\times\R^d$.
In~\cite[Sec.~7.1.3]{Burden2016-bb}, $\samp{F}$ was referred to as the \emph{sampled} vector field since it is obtained by ``sampling'' the selection functions $F_b$ that define $F$ near $\rho$,
and it was noted that the function $\samp{\phi}$ 
%(i) 
is piecewise-affine
and 
%(ii) 
it approximates the original vector field's flow $\phi$ near $\rho$.
We will leverage the algebraic properties of $\samp{\phi}$ and its relationship to $\phi$ in what follows to obtain our results.

\subsection{Time-to-impact for an {\ECr} vector field and its local approximation}
\label{sec:bg:tti}

Suppose vector field $F:\R^d\into T\R^d$ is event-selected $\Cr$ with respect to $h\in\Cr(U,\R^n)$ at $\rho\in U\subset\R^d$,
and let $\phi\in\PCr(\F,\R^d)$ be its piecewise-differentiable flow.
Then~\cite[Thm.~7]{Burden2016-bb} ensures there exists a piecewise-differentiable \emph{time-to-impact} function $\tau\in\PCr(U,\R^n)$
for which
\eqnn{\label{eq:bg:tti}
\forall x\in U, j\in\set{1,\dots,n} : \phi(\tau_j(x),x) \in H_j = h_j^{-1}(h_j(\rho)),
}
i.e.\ the point $x$ flows to the surface $H_j$ in time $\tau_j(x)$.
Similarly, applying~\cite[Thm.~7]{Burden2016-bb} to the sampled vector field $\samp{F}:\R^d\into T\R^d$ and piecewise-affine flow $\samp{\phi}:\samp{\F}\into\R^d$ associated with $F$ at $\rho$ constructed in~\secref{bg:samp} ensures there exists a piecewise-affine time-to-impact function $\samp{\tau}:\R^d\into\R^n$ for which
\eqnn{\label{eq:bg:ttisamp}
\forall x\in \R^d, j\in\set{1,\dots,n} : \samp{\phi}(\samp{\tau}_j(x),x) \in \samp{H}_j = \rho + \ker Dh_j(\rho),
}
i.e.\ the point $x$ flows to the affine subspace $\samp{H}_j$ in time $\samp{\tau}_j(x)$.

\section{Representation}
\label{sec:rep}

Our main theoretical result is an explicit representation for the Bouligand (or B-)derivative of the piecewise-differentiable flow generated by an event-selected $\Cr$ vector field.
To that end, let $F:\R^d\into T\R^d$ be an event-selected $\Cr$ vector field and $\phi:\F\into\R^d$ its piecewise-differentiable flow.
In what follows, we will assume that 
$s,t\in\R$ 
and 
$x\in \R^d$ 
are such that 
$0 < s < t$ 
and 
the vector field $F$ is $\Cr$ on $\phi([0,t]\sm\set{s},x)$.
Although a general trajectory can encounter more than one point of discontinuity for $F$, such points are isolated~\cite[Lem.~6]{Burden2016-bb}, so 
%\GC{My understanding is that this follows by the flow being $PC^r$, not just that points of non-differentiability are isolated, i.e., the chain rule holds as the flow is B-differentiable, and the B-derv follows chain rule} 
% SB agree -- trying to point out that only a finite number of function compositions are needed
the Chain Rule for B-differentiable functions~\cite[Thm.~3.1.1]{Scholtes2012-la} can be applied to 
triangulate the desired flow derivative by
composing the triangulated flow derivatives associated with each point.
Thus, 
without loss of generality,
we restrict our attention to portions of trajectories that encounter one point of discontinuity for $F$, which point lies at the intersection of $n$ surfaces of discontinuity for $F$.
We assume $n > 1$ because at least two surfaces are needed for our results to be useful:
when $n = 1$ the desired B-derivative is linear~\cite{Aizerman1958-ih}% 
%\GC{While I think this is very authoritative reference for the saltation matrix, I think van der schaff's book is a lot more readable and concrete, so maybe let's cite that as well}
% SB can you provide a specific cite so I can look at it?
, so it may be represented and employed in computations as a matrix.

The B-derivative $D_x\phi(t,x):T_x \R^d\into T_{\phi(t,x)} \R^d$ we seek is a continuous 
%\GC{I have found some readers confuse 'continuous as a map from one tangent space to another' with 'continuous as the evaluation point is varied'. Any thoughts on how to drive this point home? Maybe point out Prop 3.1.2 of Scholtes?} 
% SB I've been confused / unsure about this before as well -- I think this isn't the right place to address it, so I'll add a note in the Background
piecewise-linear function, 
so it can be parsimoniously represented using a \emph{triangulation}~\cite[Sec.~3.1]{Groff2003-pb}, that is, a combinatorial simplicial complex (as defined in~\secref{bg:pwa}) each of whose vertices are associated with a pair of (tangent) vectors -- one each in the domain and codomain of $D_x\phi(t,x)$.
We will obtain this triangulation via an indirect route:  
in~\secref{rep:tri}, we triangulate the piecewise-affine flow $\samp{\phi}$ introduced in~\secref{bg:samp};
in~\secref{rep:sampphi}, we differentiate our representation of $\samp{\phi}$ to obtain a triangulation of the B-derivative $D_x\samp{\phi}$;
in~\secref{rep:phi}, we show how the B-derivative $D_x\phi$ can be obtained from $D_x\samp{\phi}$, providing a triangulation of the desired derivative.

\subsection{Triangulation}
\label{sec:rep:tri}

\begin{figure}[t]
	\centering
	\def\svgwidth{1.3\columnwidth} 
	\resizebox{1.\columnwidth}{!}{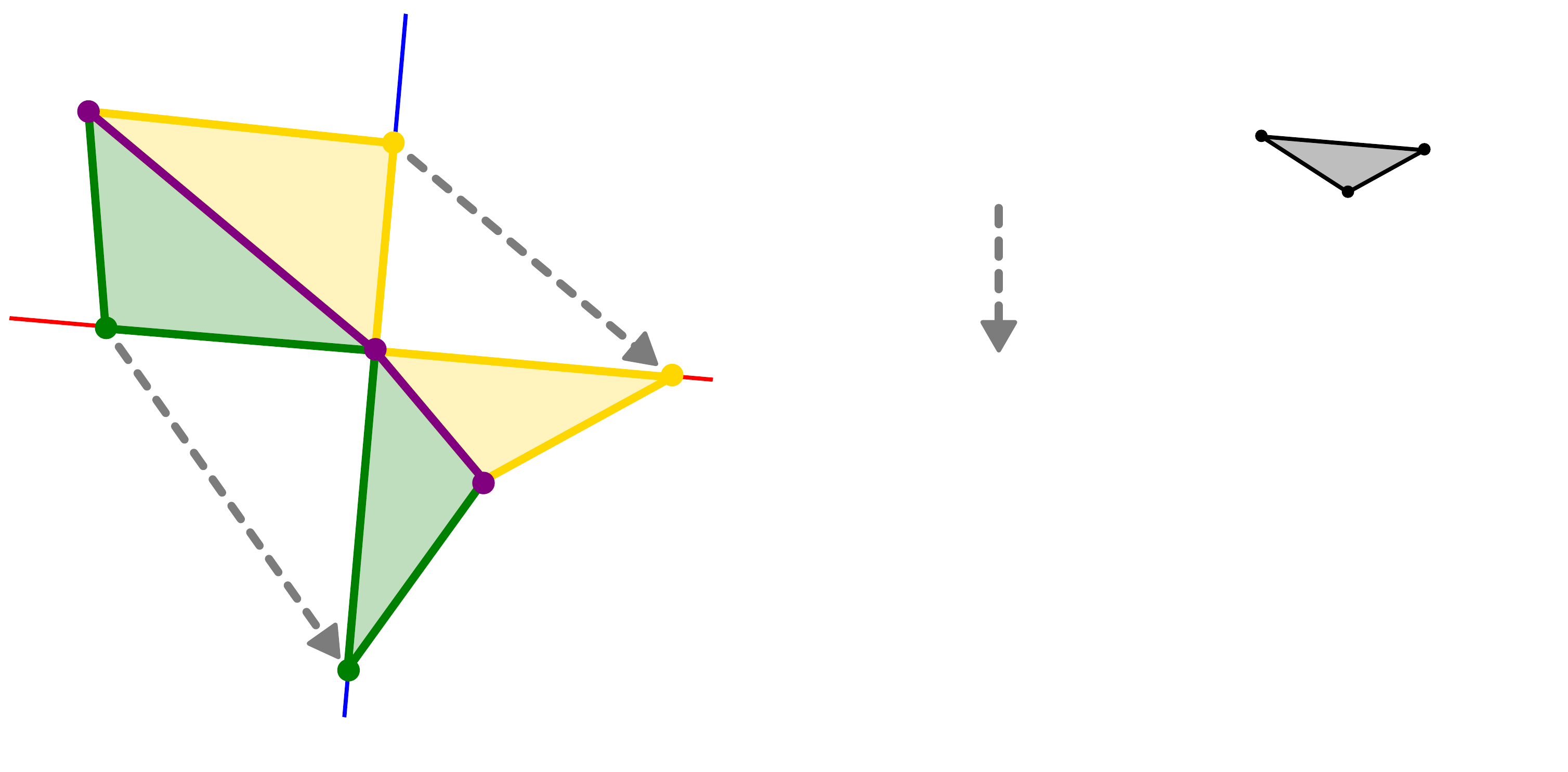}
	\caption{%
		\emph{Triangulation of the time-1 flow $\samp{\phi}_1$ of the \emph{sampled system} associated with an {\ECr} vector field.}
		(a) For each $b\in\set{-1,+1}^2$, the point $\zeta_b$ defined by~\eqref{eq:zetab} flows from $\samp{D}_{-\ones}$ to $\samp{D}_{+\ones}$ in 1 (one) unit of time via the \emph{sampled system} illustrated in~\figref{pert}(b) and defined in~\secref{bg:samp}.
		(b) The sets
		$\set{\zeta_{-\ones},\zeta_{+\ones},\zeta_{(+1,-1)}}$%
		,
		$\set{\zeta_{-\ones},\zeta_{+\ones},\zeta_{(-1,+1)}}$
		indexed by~\eqref{eq:Delta}
		define geometric simplices $\Delta_{(1,2)}^-$, $\Delta_{(2,1)}^+$
		that
		pass through subspaces $\samp{H}_1, \samp{H}_2$ in the same order.
		(c) For each $\sigma\in\set{(1,2),(2,1)}$, extending $\Delta_\sigma^-$ by direct sum with subspace $\K$ yields $\Sigma_\sigma$.
	}
	\label{fig:tri}
\end{figure}

The goal of this subsection is to triangulate the piecewise-affine flow $\samp{\phi}$ introduced in~\secref{bg:samp}.
To that end, let $\rho = \phi(s,x)$ and suppose%
\footnote{%
As observed in~\cite[Sec.~7.1.5]{Burden2016-bb}, first-order approximations of an {\ECr} vector field's {\PCr} flow are not affected by flow between surfaces that are tangent at $\rho$, so we assume such redundancy has been removed.
}
$\rank D h(\rho) = n$ so 
$%\eqnn{
\set{\delta\rho\in T_\rho \R^d : b = \sgn D h(\rho)\cdot \delta\rho}
$%} 
has nonempty interior for each $b\in\set{-1,+1}^n = \Bn$.
Letting $\K = \ker D h(\rho)\subset T_\rho \R^d$ 
denote the kernel of $D h(\rho)$ and
$\K^\perp$ 
its orthogonal complement,
for each $b\in \Bn$ there exists a unique%
\footnote{%
Here and in what follows we mildly abuse notation via the natural vector space isomorphism $\R^d \simeq T_\rho \R^d$.
}%
\footnote{%
Uniqueness is ensured by $\rank D h(\rho) = n$ since (i) $\K^\perp$ is $n$-dimensional, (ii) the rows of $D h(\rho)$ are linearly independent, and hence (iii) there are $n$ independent equations in the $n$ unknowns needed to specify $\zeta_b$ in~\eqref{eq:zetab}.
}
$\zeta_b\in\K^\perp + \set{\rho}$ such that%
%\GC{Seems like tangent and ambient space are conflated, and that a comment should be in order. The sum seems to try and address this, but to be super annoying, we don't specify then how Dh extends to this sum.}
% SB fair enough -- how's the footnote I added?
\eqnn{\label{eq:zetab}
D h_{b > 0}(\rho)(\zeta_b - \rho) = 0\ \text{and}\ D h_{b < 0}(\rho)(\zeta_b + F_b(\rho) - \rho) = 0
}
where $h_{b > 0}$ (respectively, $h_{b < 0}$) denotes the function obtained by selecting components $h_j$ of $h$ for which $b_j = +1$ (respectively, $b_j = -1$)%
%\footnote{Formally, 
%with $\card{b>0}\in\N$ denoting the number of $+1$'s in $b$,
%we define
%$h_{b > 0} = \paren{h_j}_{b_j = +1}:\R^d\into\R^{\card{b>0}}$
%(and $h_{b < 0}$ = $h_{-b > 0}$).}%
.
The vectors defined by~\eqref{eq:zetab} have special significance for the piecewise-affine flow $\samp{\phi}$ introduced in~\secref{bg:samp} (see~\figref{tri}(a)):
\eqnn{\label{eq:zetaphi}
\forall b\in\Bn : \zeta_b\in\samp{D}_{-\ones},\ \samp{\phi}(1,\zeta_b) = \zeta_b + F_b(\rho)\in\samp{D}_{+\ones}, 
}
that is, the point $\zeta_b$ lies ``before'' all event surface tangent planes and flows in $1$ (one) unit of time to $\zeta_b + F_b(\rho)$ which lies ``after'' all event surface tangent planes (neither ``before'' nor ``after'' should be interpreted strictly).
We denote the collections of these vectors as follows:
\eqnn{\label{eq:Z}
Z^- = \set{\zeta_b}_{b\in\Bn},\ Z^+ = \set{\zeta_b + F_b(\rho)}_{b\in\Bn}.
}

% fig:tri

In what follows, it will be convenient
to use an element $\sigma\in\Sn$ of the symmetric permutation group over $n$ elements 
to specify $n+1$ elements of $b\in\Bn$ as follows: 
for each $k\in\set{0,\dots,n}$, let
% SB this could use some additional explanation :)
$\sigma(\set{0,\dots,k})\subset\set{1,\dots,n}$ specify the unique $b\in\Bn$ whose $j$-th component is $+1$ if and only if $j\in\sigma(\set{0,\dots,k})$.
Note that this identification yields, with some abuse of notation, 
$\sigma(\set{0}) = -\ones$, $\sigma(\set{0,\dots,n}) = +\ones$.
Finally, note that:
%$\zeta_{ -\ones }(\rho) + F_{-\ones}(\rho)  = \rho$,
%$\zeta_{ +\ones } = \rho$,
%and we have the following two facts:
\begin{subequations}
\label{eq:zeta}
\begin{align}
\set{ \zeta_{ \sigma(\set{0,\dots,k}) } - \rho }_{k=0}^{n-1}
&\ \text{are linearly independent}; 
\label{eq:zetasigma}
\\
\set{ \zeta_{ \sigma(\set{0,\dots,k}) } + F_{ \sigma(\set{0,\dots,k}) }(\rho) - \rho }_{k=1}^{n} 
&\ \text{are linearly independent}.
\label{eq:zetasigmaFsigma}
\end{align}
\end{subequations}
The former fact~\eqref{eq:zetasigma} is easily verified in coordinates where $D h(\rho) = \left[ I_n\ 0_{n\times(d-n)} \right]$,
whence the latter fact~\eqref{eq:zetasigmaFsigma} follows from~\eqref{eq:zetasigma} and~\eqref{eq:zetaphi} via~\cite[Cor.~5(c)]{Burden2016-bb} (the time-$t$ flow of an {\ECr} vector field is a homeomorphism of the state space for all $t\in\R$).

Let $\Delta$ denote the combinatorial simplicial complex over vertex set $\Bn$ whose maximal $n$-simplices are indexed by $\sigma\in\Sn$ via
\eqnn{\label{eq:Delta}
\Delta_\sigma = \set{ \sigma(\set{0,\dots,k}) }_{k=0}^n \in \Delta
}
where we regard $\sigma(\set{0,\dots,k})$ as an element of $\Bn$ using the same abuse of notation employed in~\eqref{eq:zeta}.
By associating each vertex $b\in\Bn$ with the vector $\zeta_b\in Z^-\subset\R^d$,  
every $n$-simplex $\Delta_\sigma$ determines an $n$-dimensional geometric simplex $\Delta_\sigma^-\subset\R^d$, the dimensionality of which is ensured by~\eqref{eq:zetasigma};
similarly,~\eqref{eq:zetasigmaFsigma}
ensures that
associating each $b\in\Bn$ with $(\zeta_b + F_b(\rho))\in Z^+\subset\R^d$ determines an $n$-dimensional geometric simplex $\Delta_\sigma^+\subset\R^d$ from each $n$-simplex $\Delta_\sigma$.
Refer to~\figref{tri}(b) for an illustration when $n = 2$.
%Letting $\Delta^- = \set{\Delta_\sigma^-}_{\sigma\in\Sn}$, $\Delta^- = \set{\Delta_\sigma^-}_{\sigma\in\Sn}$ denote polyhedral subdivisions determined by the $\Delta_\sigma^\pm$'s,
%
The triple $(Z^-,Z^+,\Delta)$ parameterizes a continuous piecewise-affine homeomorphism 
%
%
%\GC{Do we need to add conv (or even carrier) here? IIRC, Groff's formula works on a ``triangulation'' -- is that equal to the carrier in our case?}
% SB yep -- that's what I realized; please check that I've fixed it
$P:\car{\Delta^-}\into\car{\Delta^+}$ using the construction from~\cite[Sec.~3.1]{Groff2003-pb} (summarized in~\secref{bg:pwa}),
where
$\car{\Delta^\pm} = \bigcup_{\sigma\in \Sn}\Delta_\sigma^\pm\subset\R^d$
denote the \emph{carriers} of the geometric simplicial complexes $\Delta^\pm$.

We now show that the piecewise-affine function $P$ constructed above is the non-linear part of 
the time-$1$ flow of the sampled system 
$\samp{\phi}_1$
restricted to $\car{\Delta^-}$. 
For each $\sigma\in\Sn$ we extend the $n$-dimensional geometric simplex $\Delta_\sigma^-$ determined by the $n$-simplex $\Delta_\sigma$ via direct sum with the $(d-n)$-dimensional subspace $\K$ to obtain a $d$-dimensional polyhedron $\Sigma_\sigma$ (see~\figref{tri}(c)),
and let $\car{\Sigma} = \bigcup_{\sigma\in\Sn}\Sigma_\sigma$. %be the resulting polyhedral subdivision; % of $\car{\Sigma} = \bigcup_{\sigma\in\Sn}\Sigma_\sigma\subset\R^d$;
% of $\K + \conv Z^-$
% \GC{carrier?} 
% SB yep -- please check I've fixed this everywhere
Note that $\K$ is a subset of the lineality space of $\Sigma_\sigma$ for each $\sigma\in\Sn$.

\begin{lemma}
  \label{lem:sampphi1}
  $\samp{\phi}_1|_{\car{\Sigma}}$ is piecewise-affine and
  \eqnn{
  %$
  \forall z\in\car{\Delta^-}, \xi\in\K : \samp{\phi}_1(z + \xi) = P(z) + \xi.
  %$
  }
\end{lemma}

\begin{proof}
  This proof will proceed in two steps:
  (i) show that $\samp{\phi}_1(z) = P(z)$ for all $z\in\car{\Delta^-}$;
  (ii) show that $\samp{\phi}_1(z + \xi) = \samp{\phi}_1(z) + \xi$ for all $z\in\car{\Delta^-}, \xi\in\K$.

  (i) Recall from~\eqref{eq:zetaphi} that $\samp{\phi}_1|_{Z^-} = P|_{Z^-}$ 
  where $Z^-$ is the vertex set for the geometric simplicial complex $\Delta^-$.
  For each $\sigma\in\Sn$ let $Z_\sigma = \set{\zeta_b}_{b\in\Delta_\sigma}$ denote the vertex set of the $n$-dimensional geometric simplex $\Delta_\sigma^-$.
  Then we claim that each $z\in\Delta_\sigma^-$ passes through the same sequence of transition surfaces as each $\zeta_b\in Z_\sigma$.
  To verify this claim, we use the piecewise-affine \emph{time-to-impact} function $\samp{\tau}:\R^d\into\R^n$ from~\secref{bg:tti}.
  Note that 
  %$\samp{\tau}(\zeta_b) = \bin{b > 0}\in\set{0,1}^n$ for each $b\in\Bn$, i.e.\ 
  $\zeta_b$ impacts affine subspace $\samp{H}_j$ at time $1$ if $b_j = -1$ and at time $0$ if $b_j = +1$, i.e.\ 
  \eqnn{\samp{\tau}_j(\zeta_b) = \pw{1, & b_j = -1; \\ 0, & b_j = +1.}}
  A convex combination $\alpha\,\zeta_b + (1-\alpha)\zeta_{b'}$, 
  $\alpha\in(0,1)$, 
  $b,b'\in\Delta_\sigma$, 
  impacts $\samp{H}_j$ at time 
  %\GC{Do not believe}
  % SB great catch!  please check that I've fixed this by restricting to b,b'\in\Delta_\sigma 
  \eqn{
    \samp{\tau}_j(\alpha\,\zeta_b + (1-\alpha)\zeta_{b'}) = 
    \pw{
      1, & b_j = -1 \wedge b'_j = -1; 
      \\
      t\in(0,1),& (b_j = +1 \wedge b'_j = -1) \vee (b_j = -1 \wedge b'_j = +1); 
      \\
      0, & b_j = +1 \wedge b'_j = +1.
    }
  }
  More generally, any point $z \in \Delta_\sigma^-$ is a convex combination of the vertices $Z_\sigma$, whence it impacts surfaces in the order prescribed by $\sigma$:
  \eqnn{\label{eq:ttisigma}
  \forall z \in \Delta_\sigma^-:
    0 \leq \samp{\tau}_{\sigma(1)}(z) 
    \leq
    \samp{\tau}_{\sigma(2)}(z) 
    \leq \cdots \leq 
    \samp{\tau}_{\sigma(n)}(z) < 1.
  }
  Thus, $\samp{\phi}_1|_{\Delta_\sigma^-}$ is affine and agrees with $P|_{\Delta_\sigma^-}$.
  Since $\car{\Delta^-} = \bigcup_{\sigma\in\Sn} \Delta_\sigma^-$, 
  we have $\samp{\phi}_1|_{\car{\Delta^-}} = P$. 

  % SB conv -> car
  (ii) We now show that the piecewise-affine map $\samp{\phi}_1$ is indifferent to $\xi\in\K = \ker Dh(\rho)$:
%  that is,
%  $\K$ is in the lineality space of $\samp{\phi}_1$ and $\samp{\phi}_1(\xi) = \xi$ for all $\xi\in\K$:
  \begin{subequations}
    \label{eq:lem:ker}
    \begin{align}
      \forall \xi\in\K, z\in\car{\Delta^-}: \samp{\phi}_1(z + \xi)
      & = \samp{\phi}_1\paren{\rho + (z + \xi - \rho)} 
      \label{eq:lem:ker:a}
      \\
      & = \samp{\phi}_1(\rho) + D\samp{\phi}_1(\rho; z + \xi - \rho) 
      \label{eq:lem:ker:b}
      \\
      & = \samp{\phi}_1(\rho) + D\samp{\phi}_1(\rho; z - \rho) + \xi 
      \label{eq:lem:ker:c}
      \\
      & = \samp{\phi}_1(z) + \xi.
      \label{eq:lem:ker:d}
    \end{align}
  \end{subequations}
  Indeed:
  \eqref{eq:lem:ker:a} since $z + \xi = \rho + (z + \xi - \rho)$;
  \eqref{eq:lem:ker:b} since $\samp{\phi}_1$ is affine on the segment $\set{\rho + \alpha\,(z+\xi-\rho) : \alpha \in [0,1]}$;
  \eqref{eq:lem:ker:c} since each piece of the continuous piecewise-linear B-derivative $D\samp{\phi}_1(\rho)$ is specified by a saltation matrix (as recalled in~\secref{bg:PCr}) that is the product of matrices of the form $(I_d + g\cdot Dh_j(\rho))$~\cite[Eqn.~(60)]{Burden2016-bb}, thus $\xi \in\K = \ker Dh(\rho)$ is transformed by $I_d$;
  \eqref{eq:lem:ker:d} for the same reason as~\eqref{eq:lem:ker:b}.
  %\GC{I think part 2 is true, and the its proof is fine, but I think we cna format it a touch differently to make it more believable.}
  % SB point taken -- revisions welcome!  I personally think (4.9b) is the biggest leap
\end{proof}

\subsection{B-derivative of $\samp{\phi}$}
\label{sec:rep:sampphi}

The goal of this subsection is to differentiate the representation of $\samp{\phi}$ from~\secref{rep:tri} to obtain a triangulation of the B-derivative 
$D\samp{\phi}_1:T_{\rho^-}\R^d \into T_{\rho^+}\R^d$
between the following two points:
\eqnn{
\rho^- = \rho-\frac{1}{2} F_{-\ones}(\rho),\ 
\rho^+ = \samp{\phi}(1,\rho^-) = \rho+\frac{1}{2} F_{+\ones}(\rho).
}

\begin{lemma}
  \label{lem:B}
  The function $B = D\samp{\phi}_1(\rho^-) : T_{\rho^-}\R^d \into T_{\rho^+}\R^d$ satisfies:
  \begin{enumerate}[leftmargin=.25in]
    \item $B$ specifies how $\samp{\phi}_1$ varies relative to $\samp{\phi}_1(\rho^-)$,
      \eqnn{\label{eq:Bsampphi}
        \forall x\in |\Sigma| : \samp{\phi}_1(x) = \samp{\phi}_1(\rho^-) + B(x - \rho^-);
      }
    \item $B$ is continuous and piecewise-linear with conical subdivision
      \eqnn{\label{eq:Sigmaprime}
        \Sigma' = \set{\Sigma'_\sigma = \cone\paren{\Sigma_\sigma - \rho^-} : \sigma\in\Sn};
      }
    \item $B|_{\Sigma'_\sigma}$ is linear for all $\sigma\in\Sn$ and 
      \eqnn{\label{eq:BSigmaprimesigma}
        \forall \delta\rho\in\Sigma'_\sigma : B(\delta\rho) = M_\sigma\cdot\delta\rho;
      }
    \item $\Li = \K + \spn{F_{-\ones}(\rho)}$
      is a 
      $(d-n+1)$-dimensional 
      lineality space for $\Sigma'$
      and
      \eqnn{\label{eq:Sigmaprimesigma}
      \forall\sigma\in\Sn : \Sigma'_\sigma = \Li + \cone\set{\Pi^\perp_{\Li}\cdot (\zeta_{\sigma(\set{0,\dots,k})}-\rho)}_{k=1}^{n-1},
      }
      where $\Pi^\perp_{\Li}$ is the orthogonal projection onto $\Li^\perp$;
    \item $B|_\Li$ is linear and 
      \eqnn{\label{eq:Bsplit}
      \forall \delta\rho\in T_{\rho^-}\R^d : B(\delta\rho) = B\paren{\Pi_\Li\cdot \delta\rho} + B\paren{\Pi^\perp_{\Li}\cdot\delta\rho},
      }
      where $\Pi_\Li$ is the orthogonal projection onto $\Li$.
  \end{enumerate}
\end{lemma}

\begin{proof}
  Each point follows from straightforward application of results in~\cite{Scholtes2012-la}:
  (1.), (2.), and (3.) are conclusions (4.), (3.), and (2.), respectively, of~\cite[Prop.~2.2.6]{Scholtes2012-la};
  (4.) follows from the definitions of lineality space~\cite[Sec.~2.1.2]{Scholtes2012-la} and the $\zeta_b$'s~\eqref{eq:zetab};
  (5.) is a restatement of~\cite[Lem.~2.3.2]{Scholtes2012-la}.
\end{proof}

\subsection{B-derivative of $\phi$}
\label{sec:rep:phi}
The goal of this subsection is to show that the piecewise-linear function $B$ triangulated in~\secref{rep:sampphi} 
%can be obtained from $D_x\samp{\phi}$.
%In the following Theorem,
%we show that the piecewise-linear function $B$ defined in Lemma~\ref{lem:B}
gives the non-linear part of the desired B-derivative $D_x \phi(t,x)$
%that the saltation matrix $M_\sigma$ from~\eqref{eq:Dxphisigma} is the matrix representation of the linear ``piece'' $B|_{\Sigma'_\sigma}$,
and%
\footnote{%
Here and in what follows we mildly abuse notation via the natural vector space isomorphisms
$\R^d \simeq T_{\rho^-} \R^d \simeq T_{\rho^+} \R^d \simeq T_\rho \R^d$.%
}
\eqnn{\label{eq:Wsigma}
W_\sigma = D_x \phi(s,x)^{-1}\paren{\Sigma'_\sigma} \subset T_x \R^d 
}
is the cone of tangent vectors where the saltation matrix 
$M_\sigma$ 
is active in~\eqref{eq:Dxphisigma}.

\begin{theorem}
  \label{thm:rep}
Suppose the vector field $F:\R^d\into T \R^d$ is event-selected $\Cr$ 
%($\ECr$) 
with respect to $h:\R^d\into\Real^n$ at $\rho$. 
Let $\phi:\F\rightarrow \R^d$ 
be the $\PCr$ flow of $F$ 
and 
$s,t\in\R$,
$x\in\R^d$
be such that 
$0 < s < t$ 
and
$F$ is $\Cr$ on $\phi( [0,t]\sm\set{s}, x)\subset \R^d$.
Then 
with
$\rho = \phi(s,x)$,
the B-derivative of the flow $\phi$ with respect to state,
$D_x\phi(t,x):T_x\R^d\into T_{\phi(t,x)}\R^d$,
%in the $\delta x\in T_x \R^d$ direction,
%$D_x\phi(t,x;\delta x)\in T_{\phi(t,x)} \R^d$,
is given by
\begin{subequations}
\label{eq:thm}
\begin{align}
  \forall \delta x \in T_x \R^d :
  D_x\phi(t,x;\delta x) 
  & = 
  D_x \phi(t-s,\rho) 
  \cdot
  B(D_x \phi(s,x)\cdot\delta x),
  \label{eq:thm:rep:B}
  \\
  \forall \delta x \in W_\sigma\subset T_x \R^d :
  D_x\phi(t,x;\delta x) 
  & = 
  D_x \phi(t-s,\rho) 
  \cdot
  M_\sigma
  \cdot
  D_x \phi(s,x)
  \cdot
  \delta x,
  \label{eq:thm:rep:Msigma}
\end{align}
\end{subequations}
where 
$B$ % = D\samp{\phi}_1(\rho^-):T_{\rho^-} \R^d\into T_\rho \R^d$ 
is the continuous piecewise-linear function from Lemma~\ref{lem:B},
$W_\sigma$ is the cone from~\eqref{eq:Wsigma},
and
$M_\sigma$ %\in\Real^{d\times d}$ 
is the saltation matrix from~\eqref{eq:Msigma}.
\end{theorem}

\begin{proof}
  Note that~\eqref{eq:thm:rep:B} follows from~\eqref{eq:thm:rep:Msigma} by~\eqref{eq:BSigmaprimesigma},
  and
  the fact that ``pieces'' of the B-derivative $D_x\phi(t,x)$ are determined by 
  the collection of saltation matrices $\set{M_\sigma}_{\sigma\in\Sn}$ was recalled in~\secref{bg:PCr}.
  Thus, to establish~\eqref{eq:thm:rep:Msigma}
  what remains to be shown is that 
  $M_\sigma$ is the active ``piece'' 
  for all $\delta x \in W_\sigma$,
  %all that remains to be shown is that 
  %$D_x\phi(t,x;\delta x) = D\phi(t-s,x) \cdot M_\sigma \cdot D\phi(s,x)\cdot \delta x$ for all
  %$\delta x \in W_\sigma = D_x\phi(s,x)^{-1}(\Sigma'_\sigma)\subset T_x\R^d$,
  i.e.\ that $\set{W_\sigma}_{\sigma\in\Sn}$ is a conical subdivision for the piecewise-linear operator $D_x\phi(t,x)$, with $W_\sigma$ as defined in~\eqref{eq:Wsigma}.  
  
  Given $\delta x\in\Int W_\sigma$ let $\delta\rho = D_x\phi(s,x)\cdot\delta x\in\Int\Sigma'_\sigma$ so that
  \eqnn{\label{eq:rep:thm:ttisamp}
    \samp{\tau}_{\sigma(1)}(\rho + \delta\rho) < 
    \samp{\tau}_{\sigma(2)}(\rho + \delta\rho) < 
    \cdots <
    \samp{\tau}_{\sigma(n)}(\rho + \delta\rho)
  }
  where $\samp{\tau}$ is the time-to-impact function for the sampled system as defined in~\eqref{eq:bg:ttisamp}.
  Note that $D_x\phi(t,x)$ is linear on $\spn F(x)$,
  %\GC{This is true, but what's the justification employed here?}
  % SB let's discuss how to explain this
  % SB this is a fundamental property of flows, but also follows immediately from the algebra
  \eqnn{
    \forall\alpha\in\R : D_x\phi(t,x; \delta x + \alpha F(x) ) = D_x\phi(t,x; \delta x) + \alpha F(\phi(t,x)),
  }
  so without loss of generality we may assume $\delta\rho \in \Int \samp{D}_{-\ones}$ by translating $\delta x$ in the $-F(x)$ direction.
  %
  % SB clarify that the threshold on \alpha depends on \delta x -- not possible to get a uniform threshold
  We claim that, for all $\alpha > 0$ sufficiently small,
  $\phi(t, x+\alpha\,\delta x)$ passes through the event surfaces with the same sequence as $\samp{\phi}(1,\rho+\alpha\,\delta\rho)$,
  i.e.\ that
  \eqnn{\label{eq:rep:thm:tti}
    \tau_{\sigma(1)}(x + \alpha\,\delta x) < 
    \tau_{\sigma(2)}(x + \alpha\,\delta x) < 
    \cdots <
    \tau_{\sigma(n)}(x + \alpha\,\delta x),
  }
  where $\tau$ is the time-to-impact function defined in~\eqref{eq:bg:tti}.
  To see this, note that
  \begin{subequations}
  \label{eq:rep:thm:Dttieq}
  \begin{align}
    \forall k\in\set{1,\dots,n} : \tau_{\sigma(k)}(x+\alpha\,\delta x) - \tau_{\sigma(k)}(x) 
    & = D\tau_{\sigma(k)}(x; \alpha\,\delta x) + \Ord{\alpha^2} 
    \label{eq:rep:thm:Dtti} \\
    & = D\samp{\tau}_{\sigma(k)}(\rho; \alpha\, \delta\rho) + \Ord{\alpha^2} 
    \label{eq:rep:thm:Dttisamp} \\
    & = \samp{\tau}_{\sigma(k)}(\rho + \alpha\,\delta\rho) - \samp{\tau}_{\sigma(k)}(\rho) + \Ord{\alpha^2}
    \label{eq:rep:thm:dttisamp} 
  \end{align}
  \end{subequations}
  where:
  \eqref{eq:rep:thm:Dtti} since $\tau$ is {\PCr};
  \eqref{eq:rep:thm:Dttisamp} since 
  $\delta\rho = D_x\phi(s,x)\cdot\delta x$
  and
  $D\tau(x;\delta x)$, $D\samp{\tau}(\rho;\delta\rho)$ are are determined by the same data, 
  namely,
  $Dh_{\sigma(k)}(\rho)$ and $F_{-\ones}(\rho)$;
  \eqref{eq:rep:thm:dttisamp} since $\delta\rho\in\Sigma'_\sigma$.
  Combining the approximation~\eqref{eq:rep:thm:Dttieq} with~\eqref{eq:rep:thm:ttisamp} yields~\eqref{eq:rep:thm:tti} as desired. 

  We conclude that $\set{W_\sigma}_{\sigma\in\Sn}$ is a conical subdivision for the piecewise-linear operator $D_x\phi(t,x)$, which verifies~\eqref{eq:thm} and completes the proof.
\end{proof}

\begin{remark}
The only non-classical part of the B-derivative of the flow in~\eqref{eq:thm:rep:B} is the piecewise-linear function $B$. 
%which 
%is the B-derivative of $\phi$ with respect to state at $(0,\rho)\in\F$,
%i.e. 
%$B(\delta\rho) = D_x\phi(0,\rho;\delta\rho)$ for all $\delta\rho\in T_\rho \R^d$.
Although there are $n!$ pieces of $B$ in general, we explicitly represent all pieces using a triangulation of $2^n$ sample points defined in~\eqref{eq:Z}, achieving a substantial reduction 
-- from factorial to ``merely'' exponential --
of the information needed to represent the first-order approximation of the flow.
Note that $B$ implicitly determines the transition sequence $\sigma$ associated with the perturbation direction $\delta x$ in~\eqref{eq:thm:rep:B},
whereas this sequence must be explicitly specified
to select the appropriate saltation matrix $M_\sigma$ in~\eqref{eq:thm:rep:Msigma}.
\end{remark}

\section{Computation}
\label{sec:comp}
We now attend to the complexity of the computational tasks required to construct or evaluate the B-derivative representation from the preceding section.
To that end,
let $F:\R^d\into T\R^d$ be an event-selected $\Cr$ vector field with respect to $h\in\Cr(\R^d,\R^n)$ and $\phi:\F\into\R^d$ its piecewise-{\Cr} flow,
and assume 
$s,t\in\R$ 
and 
$x\in \R^d$ 
are such that 
$0 < s < t$, 
$\rho = \phi(s,x)$, 
and 
the vector field $F$ is $\Cr$ on $\phi([0,t]\sm\set{s},x)$.

We seek to compute $D_x \phi(t,x;\delta x)$ given $\delta x \in T_x\R^d$.
Since~\eqref{eq:thm:rep:B} from Theorem~\ref{thm:rep} yields
\eqnn{
D_x\phi(t,x;\delta x) = D_x \phi(t-s,x) \cdot B(D_x \phi(s,x)\cdot\delta x)
}
where $B:T_\rho\R^d\into T_\rho\R^d$,
the crux of the computation is
\eqnn{\label{eq:drhop}
\delta\rho^+ = B(\delta\rho^-) %= D\samp{\phi}_1(\rho^-;\delta\rho^-)
}
where $\delta \rho^- = D_x \phi(s,x)\cdot\delta x$.
In fact, Lemma~\ref{lem:B} offers further simplification via~\eqref{eq:Bsplit}:
%\eqnn{
%B(\delta\rho^-) = B\paren{\Pi_\Li\cdot\delta\rho^-} + B\paren{\Pi^\perp_{\Li}\cdot\delta\rho^-};
%}
since $B = B\circ\Pi_\Li + B\circ\Pi_\Li^\perp$ where $B\circ\Pi_\Li$ is the linear function
\eqnn{
B\circ\Pi_\Li\cdot\delta\rho^- = \paren{I_d + \paren{F_{+\ones}(\rho) - F_{-\ones}(\rho)} \cdot \frac{F_{-\ones}(\rho)^\tr}{\norm{F_{-\ones}(\rho)}^2}}\cdot\Pi_\Li\cdot\delta\rho^-,
}
only the piecewise-linear function $B\circ\Pi_{\Li}^\perp$ 
(equivalently, the restriction $B|_{\Li^\perp}$)
requires special consideration.
In what follows, we will assume the following data, needed to construct the \emph{sampled system} illustrated in~\figref{pert}(b), is given:
linearly-independent normal vectors for the surfaces of discontinuity,
i.e.\
$Dh(\rho)\in\R^{n\times d}$ with $\rank Dh(\rho) = n$;
limiting values of the vector field at the point of intersection,
i.e.\
$F_b(\rho)\in T_\rho\R^d$ for each $b\in\Bn$;
and
F-derivatives of the continuously-differentiable parts of the flow,
i.e.\
$D_x\phi(s,x), D_x\phi(t-s,x)\in\R^{d\times d}$.

\subsection{Constructing the B-derivative}
\label{sec:comp:cons}

Lemma~\ref{lem:B} demonstrates that there are $n!$ pieces of the piecewise-linear function $B$,
namely, the collection of saltation matrices 
$\set{M_\sigma}_{\sigma\in\Sn}$ in~\eqref{eq:BSigmaprimesigma} 
that are active on the corresponding polyhedral cones in the conical subdivision 
$\Sigma' = \set{\Sigma'_\sigma}_{\sigma\in\Sn}$ in~\eqref{eq:Sigmaprime}.
These polyhedral cones are generated by the $2^{n-1}$ points $\set{\zeta_b : b\in\Bn\sm\set{-\ones,+\ones}}$ in~\eqref{eq:Sigmaprimesigma}.
For each $b\in\Bn$, the point $\zeta_b\in\K^\perp+\set{\rho}$ where $\K = \ker Dh(\rho)$ can be determined by solving the $n$ affine equations with $n$ unknowns in~\eqref{eq:zetab}.
Given $\sigma\in\Sn$, the linear piece $B|_{\Li^\perp\cap\Sigma'_\sigma}$ can be constructed using the \emph{saltation matrix}~\cite[Sec.~7.1.6]{Burden2016-bb}
since
$B(\delta\rho^-) = M_\sigma\cdot\delta\rho^-$
for all $\delta\rho^-\in\Li^\perp\cap\Sigma'_\sigma$
where%
\footnote{
We mildly abuse notation as in~\secref{rep:tri}
by using $\sigma\in\Sn$ 
to specify $n+1$ elements of $b\in\Bn$: 
for each $k\in\set{0,\dots,n}$, we let
$\sigma(\set{0,\dots,k})\subset\set{1,\dots,n}$ specify the unique $b\in\Bn$ whose $j$-th component is $+1$ if and only if $j\in\sigma(\set{0,\dots,k})$.
}
\eqnn{\label{eq:comp:Msigma}
M_\sigma = \prod_{k=0}^{n-1}\paren{I_d + \frac{\paren{ F_{\perm{k+1}}(\rho) - F_{\perm{k}}(\rho) }}{Dh_{\perm{k}}(\rho) \cdot F_{\perm{k}}(\rho)} \cdot Dh_{\perm{k}}(\rho)},
}
or using \emph{barycentric coordinates}~\cite[Eqn.~(3.1)]{Groff2003-pb} 
since 
$B(\delta\rho^-) = Z_\sigma^+\cdot(Z_\sigma^-)^\pinv\cdot\delta\rho^-$
for all $\delta\rho^-\in\Li^\perp\cap\Sigma'_\sigma$
where
\eqnn{\label{eq:comp:bary}
Z_\sigma^\pm = \mat{cccc}{z_{\sigma(\set{0,1})}^\pm & z_{\sigma(\set{0,1,2})}^\pm & \cdots & z_{\sigma(\set{0,1,\dots,n-1})}^\pm}\in\R^{d\times(n-1)},
}%
\eqnn{\label{eq:zb}
\forall b\in\Delta'_\sigma : z_b^- = \Pi^\perp_{\Li} \cdot (\zeta_b - \rho),\ z_b^+ = B|_{\Li^\perp}(z_b^-),
}%
\eqnn{\label{eq:Deltaprimesigma}
\Delta'_\sigma = \set{ \sigma(\set{0,1,\dots,k}) }_{k=1}^{n-1};
}
note that 
the pseudo-inverse $\paren{Z_\sigma^-}^\pinv$ is injective on $\Li^\perp\cap\Sigma'_\sigma$
%$\spn Z_\sigma^- = \Li^\perp$ 
by~\eqref{eq:zetasigma} and~\eqref{eq:Sigmaprimesigma}.
Although the matrices
$M_\sigma,Z_\sigma^+\cdot\paren{Z_\sigma^-}^\pinv\in\R^{d\times d}$ 
define the same linear transformation on the $(n-1)$-dimensional cone $\Li^\perp\cap\Sigma'_\sigma$,
they are generally not the same matrix.
We conclude by noting that constructing the saltation matrix in~\eqref{eq:comp:Msigma} requires \Ord{n d^2} time and \Ord{d^2} space,
whereas constructing the Barycentric coordinates in~\eqref{eq:comp:bary} requires \Ord{n^2 d^2} time and \Ord{d^2} space 
(although evaluating the expression $Z_\sigma^+\cdot\paren{Z_\sigma^-}^\pinv\cdot\delta\rho^-$ requires only \Ord{n d^2} time given $Z_\sigma^\pm$).

\subsection{Evaluating the B-derivative}
\label{sec:comp:eval}

One obvious strategy to evaluate $B$ on $\delta\rho^-\in T_\rho\R^d$
is to (i)
determine $\sigma\in\Sn$ such that $\delta\rho^-\in\Sigma'_\sigma$ 
then 
(ii) apply the corresponding saltation matrix or barycentric coordinates calculation from the preceding section.
The general formulation of (i), termed the \emph{point location} problem in the computational geometry literature, is ``essentially open''~\cite[Sec.~6.5]{De_Berg2000-kt}.
For an arrangement of $m$ hyperplanes in $\R^d$, queries can be answered in \Ord{d \log m} time at the expense of \Ord{m^d} space~\cite{Chazelle1994-jb}.
In our context, the conical subdivision $\Sigma'$ in~\eqref{eq:Sigmaprimesigma} is determined by an arrangement of $m = \Ord{n!^2}$ hyperplanes, 
so this general-purpose algorithm has time complexity $\Ord{d \log n!} = \Ord{d\, n\log n}$ and space complexity \Ord{n!^{d}}.

The relationship established by~\eqref{eq:Bsampphi} between the desired B-derivative and the flow of the \emph{sampled system} illustrated in~\figref{pert}(b) 
suggests a different strategy, summarized in~\figref{alg}, with slightly worse \Ord{n^2 d} time complexity but dramatically superior \Ord{d} space complexity.
To understand the strategy,
interpret the tangent vector $\delta\rho^-\in T_{\rho^-}\R^d$ 
as a perturbation away from the point $\rho^- = \rho-\frac{1}{2}F_{-\ones}(\rho)$ that flows through $\rho$ to $\rho^+ = \rho+\frac{1}{2} F_{+\ones}(\rho)$ in one unit of time
and observe that%
\footnote{%
This equation only holds when $\norm{\delta\rho^-}$ is small enough to ensure $\rho^-+\delta\rho^-\in\samp{D}_{-\ones}$ and $\rho^++\delta\rho^+\in\samp{D}_{+\ones}$; since the B-derivative is positively-homogeneous, we impose this restriction without loss of generality.
}
$\delta\rho^+ = \samp{\phi}_1(\rho^-+\delta\rho^-) - \rho^+ = B(\delta\rho^-)$ as in~\eqref{eq:Bsampphi}.
The flow of the sampled system $\samp{\phi}_1$ is piecewise-affine,
and can be evaluated on a given perturbation vector $\delta\rho^-$ by performing a sequence of $n$ affine projections 
(one for each of the affine subspaces $\set{\samp{H}_j}_{j=1}^n$ where $\samp{F}$ is discontinuous) 
specified by the permutation $\sigma\in\Sn$ for which $\delta\rho^-\in\Sigma'_\sigma$.
Fortuitously, the sequence $\sigma$ can be determined inductively as follows.
First, define
\eqnn{\label{eq:comp:base}
\delta t_1 &= 0,\\
\delta\rho_1 & = \delta\rho^-,\\ 
\sigma(1) &= \arg\min\set{ -\frac{Dh_j(\rho)\cdot\delta\rho_1}{Dh_j(\rho)\cdot F_{-\ones}(\rho)} : j\in\set{1,\dots,n} },\\ 
\tau_1 &= -\frac{Dh_{\sigma(1)}(\rho)\cdot\delta\rho_1}{Dh_{\sigma(1)}(\rho)\cdot F_{-\ones}(\rho)}.
}
Then for $k\in\set{1,\dots,n-1}$ inductively define
\eqnn{\label{eq:comp:ind}
\delta t_{k+1} &= \delta t_k + \tau_k,\\
\delta\rho_{k+1} &= \delta\rho_{k}+\tau_{k}\cdot F_{\sigma(\set{0,\dots,k-1})}(\rho),\\ 
\sigma(k+1) &= \arg\min\set{ -\frac{Dh_j(\rho)\cdot\delta\rho_{k+1}}{Dh_j(\rho)\cdot F_{\sigma(\set{0,\dots,k})}(\rho)} : j\in\set{1,\dots,n}\sm\sigma(\set{1,\dots,k}) },\\
\tau_{k+1} &= -\frac{Dh_{\sigma(k+1)}(\rho)\cdot\delta\rho_{k+1}}{Dh_{\sigma(k+1)}(\rho)\cdot F_{\sigma(\set{0,\dots,k})}(\rho)}.
}
Finally, set $\delta\rho^+ = \delta\rho_n - (\delta t_n+\tau_n)\cdot F_{+\ones}(\rho)$.
By construction, $\delta\rho^-\in\Sigma'_\sigma$ and $\delta\rho^+ = B(\delta\rho^-)$.
This strategy is succinctly summarized in pseudocode and sourcecode in~\figref{alg};
its time complexity is $\Ord{n^2 d}$
since there are $n$ steps in the induction and each step requires $\Ord{n}$ dot products between $d$-vectors. 
The space complexity is $\Ord{d}$
since each step in the induction requires $\Ord{d}$ storage and data from preceding steps can be forgotten or overwritten. 

We conclude by noting that, if a general-purpose algorithm is employed to solve the point location problem in \Ord{d\,n\log n} time to obtain the sequence $\sigma\in\Sn$, then the induction described in the preceding paragraph can be simplified by skipping the steps that determine $\sigma(1)$ and $\sigma(k+1)$ from~\eqref{eq:comp:base} and~\eqref{eq:comp:ind}.
This simplification reduces the time complexity of the induction to \Ord{n d}, so the overall algorithm retains the \Ord{d\,n\log n} time complexity of the general-purpose point-location algorithm (at the expense of the superexponential \Ord{n!^d} space complexity of the point location algorithm).
We are pessimistic these asymptotic complexities can be improved in general.

\section{Conclusion}
\label{sec:conc}

We constructed a representation for the \emph{Bouligand} (or \emph{B-})derivative of the \emph{piecewise-{\Cr}} ({\PCr}) \emph{flow} generated by an \emph{event-selected {\Cr}} ({\ECr}) \emph{vector field} and applied the representation to derive a polynomial-time algorithm to evaluate the B-derivative on a given tangent vector.
Our results provide a foundation that may support future work generalizing classical analysis and synthesis techniques for smooth control systems to the class of nonsmooth systems considered here.
In particular, we envision applying our results to design and control the class of mechanical systems subject to unilateral constraints that arise in models of robot locomotion and manipulation.


\begin{thebibliography}{10}

\bibitem{Aguilar2015-fy}
{\sc J.~Aguilar and D.~I. Goldman}, {\em Robophysical study of jumping dynamics
  on granular media}, Nature physics, 12 (2015), p.~nphys3568,
  \url{https://doi.org/10.1038/nphys3568}.

\bibitem{Aizerman1958-ih}
{\sc M.~A. Aizerman and F.~R. Gantmacher}, {\em Determination of stability by
  linear approximation of a periodic solution of a system of differential
  equations with discontinuous right-hand sides}, The Quarterly Journal of
  Mechanics and Applied Mathematics, 11 (1958), pp.~385--398,
  \url{https://doi.org/10.1093/qjmam/11.4.385}.

\bibitem{Ballard2000-ui}
{\sc P.~Ballard}, {\em The dynamics of discrete mechanical systems with perfect
  unilateral constraints}, Archive for Rational Mechanics and Analysis, 154
  (2000), pp.~199--274, \url{https://doi.org/10.1007/s002050000105}.

\bibitem{Bertsekas1999-an}
{\sc D.~P. Bertsekas}, {\em Nonlinear Programming}, Athena Scientific, 2nd~ed.,
  1999.

\bibitem{Bizzarri2013-jh}
{\sc F.~Bizzarri, A.~Brambilla, and G.~Storti~Gajani}, {\em Lyapunov exponents
  computation for hybrid neurons}, Journal of Computational Neuroscience, 35
  (2013), pp.~201--212, \url{https://doi.org/10.1007/s10827-013-0448-6}.

\bibitem{Burden2016-bb}
{\sc S.~A. Burden, S.~S. Sastry, D.~E. Koditschek, and S.~Revzen}, {\em
  Event-selected vector field discontinuities yield piecewise-differentiable
  flows}, SIAM Journal on Applied Dynamical Systems, 15 (2016), pp.~1227--1267,
  \url{https://doi.org/10.1137/15M1016588}.

\bibitem{Chazelle1994-jb}
{\sc B.~Chazelle and J.~Friedman}, {\em Point location among hyperplanes and
  unidirectional ray-shooting}, Computational Geometry, 4 (1994), pp.~53--62,
  \url{https://doi.org/10.1016/0925-7721(94)90009-4}.

\bibitem{Collins2015-sk}
{\sc S.~H. Collins, M.~B. Wiggin, and G.~S. Sawicki}, {\em Reducing the energy
  cost of human walking using an unpowered exoskeleton}, Nature, 522 (2015),
  pp.~212--215, \url{https://doi.org/10.1038/nature14288}.

\bibitem{De_Berg2000-kt}
{\sc M.~de~Berg, M.~van Kreveld, M.~Overmars, and O.~Schwarzkopf}, {\em
  Computational Geometry: Algorithms and Applications}, Springer, 2000.

\bibitem{Bernardo2008piecewise}
{\sc M.~di~Bernardo, C.~Budd, A.~R. Champneys, and P.~Kowalczyk}, {\em
  Piecewise-smooth dynamical systems: theory and applications}, vol.~163,
  Springer Science $+$ Business Media, 2008.

\bibitem{Dieci2011-ps}
{\sc L.~Dieci and L.~Lopez}, {\em Fundamental matrix solutions of piecewise
  smooth differential systems}, Mathematics and Computers in Simulation, 81
  (2011), pp.~932--953, \url{https://doi.org/10.1016/j.matcom.2010.10.012}.

\bibitem{Elandt2019-yz}
{\sc R.~Elandt, E.~Drumwright, M.~Sherman, and A.~Ruina}, {\em A pressure field
  model for fast, robust approximation of net contact force and moment between
  nominally rigid objects}, in {IEEE/RSJ} International Conference on
  Intelligent Robots and Systems ({IROS}), Nov. 2019, pp.~8238--8245,
  \url{https://doi.org/10.1109/IROS40897.2019.8968548}.

\bibitem{Eldering2016-oj}
{\sc J.~Eldering and H.~Jacobs}, {\em The role of symmetry and dissipation in
  biolocomotion}, SIAM Journal on Applied Dynamical Systems, 15 (2016),
  pp.~24--59, \url{https://doi.org/10.1137/140970914}.

\bibitem{Filippov1988-nh}
{\sc A.~F. Filippov}, {\em Differential equations with discontinuous righthand
  sides}, Springer, 1988.

\bibitem{Groff2003-pb}
{\sc R.~E. Groff}, {\em Piecewise linear homeomorphisms for approximation of
  invertible maps}, PhD thesis, University of Michigan, 2003.

\bibitem{Hatcher2002-ec}
{\sc A.~Hatcher}, {\em {Algebraic topology}}, Cambridge University Press, 2002.

\bibitem{Hespanha2009-nf}
{\sc J.~P. Hespanha}, {\em Linear systems theory}, Princeton University Press,
  2009.

\bibitem{Hiskens2000-ps}
{\sc I.~A. Hiskens and M.~A. Pai}, {\em Trajectory sensitivity analysis of
  hybrid systems}, IEEE Transactions on Circuits and Systems I: Fundamental
  Theory and Applications, 47 (2000), pp.~204--220,
  \url{https://doi.org/10.1109/81.828574}.

\bibitem{Ivanov1998-ff}
{\sc A.~P. Ivanov}, {\em The stability of periodic solutions of discontinuous
  systems that intersect several surfaces of discontinuity}, Journal of Applied
  Mathematics and Mechanics, 62 (1998), pp.~677--685,
  \url{https://doi.org/10.1016/S0021-8928(98)00087-2}.

\bibitem{Jeffrey2014-nt}
{\sc M.~R. Jeffrey}, {\em Dynamics at a switching intersection: hierarchy,
  isonomy, and multiple sliding}, SIAM Journal on Applied Dynamical Systems, 13
  (2014), pp.~1082--1105, \url{https://doi.org/10.1137/13093368X}.

\bibitem{Johnson2016-nh}
{\sc A.~M. Johnson, S.~A. Burden, and D.~E. Koditschek}, {\em A hybrid systems
  model for simple manipulation and self-manipulation systems}, The
  International Journal of Robotics Research, 35 (2016), pp.~1354--1392,
  \url{https://doi.org/10.1177/0278364916639380}.

\bibitem{Leine2013dynamics}
{\sc R.~I. Leine and H.~Nijmeijer}, {\em Dynamics and bifurcations of
  non-smooth mechanical systems}, vol.~18, Springer Science $+$ Business Media,
  2013.

\bibitem{Ljung1999-ee}
{\sc L.~Ljung}, {\em System identification: theory for the user},
  Prentice-Hall, 1999.

\bibitem{Lotstedt1982-ea}
{\sc P.~L{\"{o}}tstedt}, {\em Mechanical systems of rigid bodies subject to
  unilateral constraints}, SIAM Journal on Applied Mathematics, 42 (1982),
  pp.~281--296, \url{https://doi.org/10.1137/0142022}.

\bibitem{Oliphant2006-ar}
{\sc T.~E. Oliphant}, {\em A guide to NumPy}, vol.~1, Trelgol Publishing USA,
  2006.

\bibitem{Pace2017-tt}
{\sc A.~M. Pace and S.~A. Burden}, {\em Piecewise-differentiable trajectory
  outcomes in mechanical systems subject to unilateral constraints}, in Hybrid
  Systems: Computation and Control ({HSCC}), 2017, pp.~243--252,
  \url{https://doi.org/10.1145/3049797.3049807}.

\bibitem{Parker1989-df}
{\sc T.~S. Parker and L.~O. Chua}, {\em Practical numerical algorithms for
  chaotic systems}, Springer, 1989.

\bibitem{Polak1997-xd}
{\sc E.~Polak}, {\em Optimization: algorithms and consistent approximations},
  Springer-Verlag, 1997.

\bibitem{Pontryagin1962-rn}
{\sc L.~S. Pontryagin, V.~G. Boltyanskii, R.~V. Gamkrelidze, and E.~F.
  Mishchenko}, {\em The mathematical theory of optimal processes (translated by
  {KN} Trirogoff)}, John Wiley \& Sons, 1962.

\bibitem{Python_Software_Foundation_undated-lp}
{\sc {Python Software Foundation}}, {\em Python language reference, version
  3.7}, \url{https://docs.python.org/release/3.7.0/}.

\bibitem{Ralph1997-yr}
{\sc D.~Ralph and S.~Scholtes}, {\em Sensitivity analysis of composite
  piecewise smooth equations}, Mathematical Programming, 76 (1997),
  pp.~593--612, \url{https://doi.org/10.1007/BF02614400}.

\bibitem{Robinson1987-va}
{\sc S.~M. Robinson}, {\em Local structure of feasible sets in nonlinear
  programming, part {III}: stability and sensitivity}, Nonlinear Analysis and
  Optimization, 30 (1987), pp.~45--66,
  \url{https://doi.org/10.1007/BFb0121154}.

\bibitem{Rockafellar2003-nd}
{\sc R.~T. Rockafellar}, {\em A property of piecewise smooth functions},
  Computational Optimization and Applications, 25 (2003), pp.~247--250,
  \url{https://doi.org/10.1023/A:1022921624832}.

\bibitem{Sastry1989-zr}
{\sc S.~Sastry and M.~Bodson}, {\em {Adaptive control: stability, convergence,
  and robustness}}, Prentice Hall, 1989.

\bibitem{Sastry1999-ei}
{\sc S.~S. Sastry}, {\em Nonlinear Systems: Analysis, Stability, and Control},
  Springer, 1999.

\bibitem{Scholtes2012-la}
{\sc S.~Scholtes}, {\em Introduction to piecewise differentiable equations},
  Springer-Verlag, 2012, \url{https://doi.org/10.1007/978-1-4614-4340-7}.

\bibitem{Simic2005-fv}
{\sc S.~N. Simic, K.~H. Johansson, J.~Lygeros, and S.~Sastry}, {\em {Towards a
  geometric theory of hybrid systems}}, Dynamics of Continuous, Discrete \&
  Impulsive Systems. Series B. Applications \& Algorithms, 12 (2005),
  pp.~649--687,
  \url{http://www.diva-portal.org/smash/record.jsf?pid=diva2:437188}.

\bibitem{Sreenath2011-lb}
{\sc K.~Sreenath, H.-W. Park, I.~Poulakakis, and J.~W. Grizzle}, {\em {A
  Compliant Hybrid Zero Dynamics Controller for Stable, Efficient and Fast
  Bipedal Walking on {MABEL}}}, The International Journal of Robotics Research,
  30 (2011), pp.~1170--1193, \url{https://doi.org/10.1177/0278364910379882}.

\bibitem{Tornambe1999-dr}
{\sc A.~Tornambe}, {\em Modeling and control of impact in mechanical systems:
  theory and experimental results}, IEEE Transactions on Automatic Control, 44
  (1999), pp.~294--309, \url{https://doi.org/10.1109/9.746255}.

\bibitem{Utkin1977-rm}
{\sc V.~Utkin}, {\em Variable structure systems with sliding modes}, IEEE
  Transactions on Automatic Control, 22 (1977), pp.~212--222,
  \url{https://doi.org/10.1109/TAC.1977.1101446}.

\end{thebibliography}
\end{document}